\documentclass[10.5pt]{article}
\usepackage{amsmath}
\usepackage{amssymb}
\usepackage{amscd}
\usepackage{xspace}
\usepackage{verbatim}
\newcommand{\mysection}[1]{
\section{#1}\setcounter{equation}{0}}
\title{\bf Reduced measures associated to parabolic problems}
\author{{\bf Waad Al Sayed}\\
{\small Department of Mathematics, Universit\'e Fran\c{c}ois Rabelais,  Tours,  FRANCE}\\[2mm]
{\bf Mustapha Jazar}\\
{\small Department of Mathematics, Universit\'e Libanaise, Beyrouth, LIBAN}\\[2mm]
 {\bf Laurent V\'eron}\\
{\small Department of Mathematics, Universit\'e Fran\c{c}ois Rabelais,  Tours,  FRANCE}
}

\date{}
\begin{document}
\maketitle
{\small {\bf Abstract} We study the existence and the properties of the reduced measures for the parabolic  equations  $\partial_tu-\Delta u+g(u)=0$ in $\Omega\times (0,\infty)$ subject to the conditions
($P$): $u=0$ on $\partial\Omega\times (0,\infty)$, $u(x,0)=\mu$ and ($P'$): $u=\mu'$ on $\partial\Omega\times (0,\infty)$, $u(x,0)=0$ where $\mu$ and $\mu'$ are positive Radon measures and $g$ a continuous nondecreasing function.
}

\vspace{1mm}
\hspace{.05in}

\newcommand{\txt}[1]{\;\text{ #1 }\;}
\newcommand{\tbf}{\textbf}
\newcommand{\tit}{\textit}
\newcommand{\tsc}{\textsc}
\newcommand{\trm}{\textrm}
\newcommand{\mbf}{\mathbf}
\newcommand{\mrm}{\mathrm}
\newcommand{\bsym}{\boldsymbol}
\newcommand{\scs}{\scriptstyle}
\newcommand{\sss}{\scriptscriptstyle}
\newcommand{\txts}{\textstyle}
\newcommand{\dsps}{\displaystyle}
\newcommand{\fnz}{\footnotesize}
\newcommand{\scz}{\scriptsize}
\newcommand{\be}{
\begin{equation}
}
\newcommand{\bel}[1]{
\begin{equation}
\label{#1}}
\newcommand{\ee}{
\end{equation}
}
\newcommand{\eqnl}[2]{
\begin{equation}
\label{#1}{#2}
\end{equation}
}
\newtheorem{subn}{\name}
\renewcommand{\thesubn}{}
\newcommand{\bsn}[1]{\def\name{#1}
\begin{subn}}
\newcommand{\esn}{
\end{subn}}
\newtheorem{sub}{\name}[section]
\newcommand{\dn}[1]{\def\name{#1}}   
\newcommand{\bs}{
\begin{sub}}
\newcommand{\es}{
\end{sub}}
\newcommand{\bsl}[1]{
\begin{sub}\label{#1}}
\newcommand{\bth}[1]{\def\name{Theorem}
\begin{sub}\label{t:#1}}
\newcommand{\blemma}[1]{\def\name{Lemma}
\begin{sub}\label{l:#1}}
\newcommand{\bcor}[1]{\def\name{Corollary}
\begin{sub}\label{c:#1}}
\newcommand{\bdef}[1]{\def\name{Definition}
\begin{sub}\label{d:#1}}
\newcommand{\bprop}[1]{\def\name{Proposition}
\begin{sub}\label{p:#1}}
\newcommand{\R}{\eqref}
\newcommand{\rth}[1]{Theorem~\ref{t:#1}}
\newcommand{\rlemma}[1]{Lemma~\ref{l:#1}}
\newcommand{\rcor}[1]{Corollary~\ref{c:#1}}
\newcommand{\rdef}[1]{Definition~\ref{d:#1}}
\newcommand{\rprop}[1]{Proposition~\ref{p:#1}}
\newcommand{\BA}{
\begin{array}}
\newcommand{\EA}{
\end{array}}
\newcommand{\BAN}{\renewcommand{\arraystretch}{1.2}
\setlength{\arraycolsep}{2pt}
\begin{array}}
\newcommand{\BAV}[2]{\renewcommand{\arraystretch}{#1}
\setlength{\arraycolsep}{#2}
\begin{array}}
\newcommand{\BSA}{
\begin{subarray}}
\newcommand{\ESA}{
\end{subarray}}
\newcommand{\BAL}{
\begin{aligned}}
\newcommand{\EAL}{
\end{aligned}}
\newcommand{\BALG}{
\begin{alignat}}
\newcommand{\EALG}{
\end{alignat}}
\newcommand{\BALGN}{
\begin{alignat*}}
\newcommand{\EALGN}{
\end{alignat*}}
\newcommand{\note}[1]{\textit{#1.}\hspace{2mm}}
\newcommand{\Proof}{\note{Proof}}
\newcommand{\qeda}{\hspace{10mm}\hfill $\square$}
\newcommand{\qed}{\\
${}$ \hfill $\square$}
\newcommand{\Remark}{\note{Remark}}
\newcommand{\modin}{$\,$\\
[-4mm] \indent}
\newcommand{\forevery}{\quad \forall}
\newcommand{\set}[1]{\{#1\}}
\newcommand{\setdef}[2]{\{\,#1:\,#2\,\}}
\newcommand{\setm}[2]{\{\,#1\mid #2\,\}}
\newcommand{\lra}{\longrightarrow}
\newcommand{\lla}{\longleftarrow}
\newcommand{\llra}{\longleftrightarrow}
\newcommand{\Lra}{\Longrightarrow}
\newcommand{\Lla}{\Longleftarrow}
\newcommand{\Llra}{\Longleftrightarrow}
\newcommand{\warrow}{\rightharpoonup}
\newcommand{
\paran}[1]{\left (#1 \right )}
\newcommand{\sqbr}[1]{\left [#1 \right ]}
\newcommand{\curlybr}[1]{\left \{#1 \right \}}
\newcommand{\abs}[1]{\left |#1\right |}
\newcommand{\norm}[1]{\left \|#1\right \|}
\newcommand{
\paranb}[1]{\big (#1 \big )}
\newcommand{\lsqbrb}[1]{\big [#1 \big ]}
\newcommand{\lcurlybrb}[1]{\big \{#1 \big \}}
\newcommand{\absb}[1]{\big |#1\big |}
\newcommand{\normb}[1]{\big \|#1\big \|}
\newcommand{
\paranB}[1]{\Big (#1 \Big )}
\newcommand{\absB}[1]{\Big |#1\Big |}
\newcommand{\normB}[1]{\Big \|#1\Big \|}

\newcommand{\thkl}{\rule[-.5mm]{.3mm}{3mm}}
\newcommand{\thknorm}[1]{\thkl #1 \thkl\,}
\newcommand{\trinorm}[1]{|\!|\!| #1 |\!|\!|\,}
\newcommand{\bang}[1]{\langle #1 \rangle}
\def\angb<#1>{\langle #1 \rangle}
\newcommand{\vstrut}[1]{\rule{0mm}{#1}}
\newcommand{\rec}[1]{\frac{1}{#1}}
\newcommand{\opname}[1]{\mbox{\rm #1}\,}
\newcommand{\supp}{\opname{supp}}
\newcommand{\dist}{\opname{dist}}
\newcommand{\myfrac}[2]{{\displaystyle \frac{#1}{#2} }}
\newcommand{\myint}[2]{{\displaystyle \int_{#1}^{#2}}}
\newcommand{\mysum}[2]{{\displaystyle \sum_{#1}^{#2}}}
\newcommand {\dint}{{\displaystyle \int\!\!\int}}
\newcommand{\q}{\quad}
\newcommand{\qq}{\qquad}
\newcommand{\hsp}[1]{\hspace{#1mm}}
\newcommand{\vsp}[1]{\vspace{#1mm}}
\newcommand{\ity}{\infty}
\newcommand{\prt}{\partial}
\newcommand{\sms}{\setminus}
\newcommand{\ems}{\emptyset}
\newcommand{\ti}{\times}
\newcommand{\pr}{^\prime}
\newcommand{\ppr}{^{\prime\prime}}
\newcommand{\tl}{\tilde}
\newcommand{\sbs}{\subset}
\newcommand{\sbeq}{\subseteq}
\newcommand{\nind}{\noindent}
\newcommand{\ind}{\indent}
\newcommand{\ovl}{\overline}
\newcommand{\unl}{\underline}
\newcommand{\nin}{\not\in}
\newcommand{\pfrac}[2]{\genfrac{(}{)}{}{}{#1}{#2}}

\def\ga{\alpha}     \def\gb{\beta}       \def\gg{\gamma}
\def\gc{\chi}       \def\gd{\delta}      \def\ge{\epsilon}
\def\gth{\theta}                         \def\vge{\varepsilon}
\def\gf{\phi}       \def\vgf{\varphi}    \def\gh{\eta}
\def\gi{\iota}      \def\gk{\kappa}      \def\gl{\lambda}
\def\gm{\mu}        \def\gn{\nu}         \def\gp{\pi}
\def\vgp{\varpi}    \def\gr{\rho}        \def\vgr{\varrho}
\def\gs{\sigma}     \def\vgs{\varsigma}  \def\gt{\tau}
\def\gu{\upsilon}   \def\gv{\vartheta}   \def\gw{\omega}
\def\gx{\xi}        \def\gy{\psi}        \def\gz{\zeta}
\def\Gg{\Gamma}     \def\Gd{\Delta}      \def\Gf{\Phi}
\def\Gth{\Theta}
\def\Gl{\Lambda}    \def\Gs{\Sigma}      \def\Gp{\Pi}
\def\Gw{\Omega}     \def\Gx{\Xi}         \def\Gy{\Psi}

\def\CS{{\mathcal S}}   \def\CM{{\mathcal M}}   \def\CN{{\mathcal N}}
\def\CR{{\mathcal R}}   \def\CO{{\mathcal O}}   \def\CP{{\mathcal P}}
\def\CA{{\mathcal A}}   \def\CB{{\mathcal B}}   \def\CC{{\mathcal C}}
\def\CD{{\mathcal D}}   \def\CE{{\mathcal E}}   \def\CF{{\mathcal F}}
\def\CG{{\mathcal G}}   \def\CH{{\mathcal H}}   \def\CI{{\mathcal I}}
\def\CJ{{\mathcal J}}   \def\CK{{\mathcal K}}   \def\CL{{\mathcal L}}
\def\CT{{\mathcal T}}   \def\CU{{\mathcal U}}   \def\CV{{\mathcal V}}
\def\CZ{{\mathcal Z}}   \def\CX{{\mathcal X}}   \def\CY{{\mathcal Y}}
\def\CW{{\mathcal W}} \def\CQ{{\mathcal Q}}
\def\BBA {\mathbb A}   \def\BBb {\mathbb B}    \def\BBC {\mathbb C}
\def\BBD {\mathbb D}   \def\BBE {\mathbb E}    \def\BBF {\mathbb F}
\def\BBG {\mathbb G}   \def\BBH {\mathbb H}    \def\BBI {\mathbb I}
\def\BBJ {\mathbb J}   \def\BBK {\mathbb K}    \def\BBL {\mathbb L}
\def\BBM {\mathbb M}   \def\BBN {\mathbb N}    \def\BBO {\mathbb O}
\def\BBP {\mathbb P}   \def\BBR {\mathbb R}    \def\BBS {\mathbb S}
\def\BBT {\mathbb T}   \def\BBU {\mathbb U}    \def\BBV {\mathbb V}
\def\BBW {\mathbb W}   \def\BBX {\mathbb X}    \def\BBY {\mathbb Y}
\def\BBZ {\mathbb Z}

\def\GTA {\mathfrak A}   \def\GTB {\mathfrak B}    \def\GTC {\mathfrak C}
\def\GTD {\mathfrak D}   \def\GTE {\mathfrak E}    \def\GTF {\mathfrak F}
\def\GTG {\mathfrak G}   \def\GTH {\mathfrak H}    \def\GTI {\mathfrak I}
\def\GTJ {\mathfrak J}   \def\GTK {\mathfrak K}    \def\GTL {\mathfrak L}
\def\GTM {\mathfrak M}   \def\GTN {\mathfrak N}    \def\GTO {\mathfrak O}
\def\GTP {\mathfrak P}   \def\GTR {\mathfrak R}    \def\GTS {\mathfrak S}
\def\GTT {\mathfrak T}   \def\GTU {\mathfrak U}    \def\GTV {\mathfrak V}
\def\GTW {\mathfrak W}   \def\GTX {\mathfrak X}    \def\GTY {\mathfrak Y}
\def\GTZ {\mathfrak Z}   \def\GTQ {\mathfrak Q}

\font\Sym= msam10 
\def\SYM#1{\hbox{\Sym #1}}
\newcommand{\bdw}{\prt\Gw\xspace}
\medskip
\mysection {Introduction}

\setcounter{equation}{0}
Let $\Gw$ be a bounded domain of $\BBR^N$, $N\geq 1$ and $g$ a nondecreasing continuous function defined on $\BBR$ and vanishing on $(-\infty,0]$. This article is concerned with the following question: {\it Given a positive Radon measure $\gn$ on
$\Gw$, does it exist a largest Radon measure $\gm$ below it for which the initial value problem
\begin {equation}\label {Eq1-1}\left\{\BA {l}
\prt_tu-\Gd u+g(u)=0\quad \text {in } Q_T:=\Gw\ti(0,T)\\[2mm]
\phantom{\prt_t-\Gd u+g(u)}
u=0\quad \text {in } \prt_\ell Q_T:=\prt\Gw\ti(0,T)\\[2mm]
\phantom{\Gd u+g(u)}
u(.,0)=\gm\quad \text {in } \Gw.
\EA\right.
\end {equation}
admits a solution}?  Whenever $\gm$ exists, it is called the {\it reduced measure} associated to $\gn$. A positive Radon measure for which (\ref{Eq1-1}) is solvable is called a  {\it good measure}. This type of problems is now well understood for nonlinear elliptic equations. This relaxation phenomenon appeared in the measure framework in the paper \cite {Va}  by Vazquez dealing with solving the problem
\begin {equation}\label {Eq1-2}
-\Gd u+e^{au}=\gm\quad \text {in }\BBR^2.
\end {equation}
He proved that the reduced measures is the sum of the non-atomic part of $\gm$ and the
atomic part where the coefficients of the Dirac masses at any atom $a$ are truncated from above at the value $2\gp/a$. Recently the general relaxation problems for the nonlinear elliptic  equations
\begin {equation}\label {Eq1-3}\left\{\BA {l}
-\Gd u+g(u)=\gm\quad \text {in }\Gw\subset\BBR^N\\[2mm]
\phantom{-\Gd u+g(u}u=0\quad \text {in }\prt\Gw
\EA\right.
\end {equation}
and
\begin {equation}\label {Eq1-4}\left\{\BA {l}
-\Gd u+g(u)=0\quad \text {in }\Gw\subset\BBR^N\\[2mm]
\phantom{-\Gd u+g(u}u=\gm\quad \text {in }\prt\Gw
\EA\right.
\end {equation}
are studied respectively by Brezis, Marcus and Ponce \cite{BMP} and Brezis and Ponce
\cite {BP}. They prove the existence of a reduced measure $\gm^*$ and study its  properties, in particular its continuity properties with respect to the capacity $W^{1,2}$ for
problem (\ref {Eq1-3}), or the (N-1)-dimensional Hausdorff measure for
problem (\ref {Eq1-4}). \smallskip

In this article we study the initial value problem in this perspective and we prove that for any  positive bounded Radon measure $\gm$ in $\Gw$ there exists a largest measure $\gm^*$, smaller than $\gm$ such that (\ref{Eq1-1}) is solvable. We study the set of good measures relative to $g$ and prove that any good measure is absolutely continuous with respect to the Hausdorff measure $H^N$. In a similar way we study the Cauchy-Dirichlet problem
\begin {equation}\label {Eq1-5}\left\{\BA {l}
\prt_tu-\Gd u+g(u)=0\quad \text {in } Q_T:=\Gw\ti(0,T)\\[2mm]
\phantom{\prt_t-\Gd u+g(u)}
u=\gm\quad \text {in } \prt_\ell Q_T:=\prt\Gw\ti(0,T)\\[2mm]
\phantom{\Gd u+g(u)}
u(.,0)=0\quad \text {in } \Gw,
\EA\right.
\end {equation}
and we prove that the reduced measure is absolutely continuous with respect to the same
Hausdorff measure $H^N$.

The proof of many results here follows the ideas borrowed from  the theory of reduced measures for elliptic equations as it is developed in  \cite{BMP} and \cite {BP}. We choose to expose them for the sake of completeness.
\section{Initial value problem}
\setcounter{equation}{0}
In this section $\Gw$ is a bounded domain in $\BBR^N$ and $\gr(x)=\dist (x,\prt\Gw)$. We denote by $\mathfrak M(\Gw)$  the set of Radon measures in $\Gw$ and, for $\ga\in\BBR$, by $\mathfrak M^\ga(\Gw)$ the subset of $\gm\in\frak M(\Gw)$ satisfying
$$\myint{\Gw}{}\gr^\ga(x)\,d\abs\gm<\infty.
$$
Thus $\mathfrak M^\ga_+(\Gw)$ is the positive cone and $\mathfrak M^0_+(\Gw)$ the set of bounded measures. For $q\in [1,\infty)$, we denote by
$L^q_{\gr^\ga}(\Gw)$ the corresponding weighted Lebesgue spaces.
For $0\leq \gt<\gs\leq T$ we set
$Q_{\gt,\gs}:=\Gw\ti (\gt,\gs)$, $Q_{\gs}:=\Gw\ti (0,\gs)$ and denote by
$\prt_\ell Q_{\gt,\gs}:=\prt\Gw\ti (\gt,\gs]$ and  $\prt_\ell Q_{\gs}:=\prt\Gw\ti (0,\gs]$ the lateral boundary of these sets. Throughout this paper we make the following assumption on $g$
\begin{equation}\label{g1}
g \text{ is a nondecreasing continuous function defined
on } \BBR \text{ and vanishing on }(-\infty, 0].
\end{equation}
\bdef {ws} Let $\gm\in \mathfrak M_+^1(\Gw)$. A function $u\in L^1(Q_T)$ is a weak solution of
(\ref{Eq1-1}) in $Q_T$ if $g(u)\in L_\gr^1(Q_T)$ and
\begin {equation}\label {Eq2-1}
\dint_{Q_T}\left(-u\prt_t\gz-u\Gd\gz+\gz\,g(u) \right)dx=\myint{\Gw}{}\gz\,d\gm,
\end {equation}
for all $\gz\in C_{\ell,0}^{2,1}(\overline Q_T)$, which is the space of functions in $C^{2,1}(\overline Q_T)$ which vanish on
$\prt\Gw\ti [0,T]\cup \overline\Gw\ti\{T\}$.
\es
We define in a similar way a weak subsolution (resp. supersolution) of (\ref{Eq1-1}) by imposing
the same integrability conditions on $u$ and $g(u)$ and
\begin {equation}\label {Eq2-1'}
\dint_{Q_T}\left(-u\prt_t\gz-u\Gd\gz+\gz\,g(u) \right)dx\,dt\leq\myint{\Gw}{}\gz\,d\gm,
\end {equation}
resp.
\begin {equation}\label {Eq2-1''}
\dint_{Q_T}\left(-u\prt_t\gz-u\Gd\gz+\gz\,g(u) \right)dx\,dt\geq\myint{\Gw}{}\gz\,d\gm,
\end {equation}
for all positive test functions in the same space. More generally we define a subsolution (resp. supersolution) of equation
\begin{equation}\label{equ}
\prt_tu-\Gd u+g(u)=0\quad\text{in } Q_{T}
\end {equation}
as a function $u\in L_{loc}^1(Q_{T})$ such that $g(u)\in L_{loc}^1(Q_{T})$ and
\begin {equation}\label {Eq2-1'1}
\dint_{Q_T}\left(-u\prt_t\gz-u\Gd\gz+\gz\,g(u) \right)dx\,dt\leq 0,
\end {equation}
resp.
\begin {equation}\label {Eq2-1''1}
\dint_{Q_T}\left(-u\prt_t\gz-u\Gd\gz+\gz\,g(u) \right)dx\,dt\geq 0,
\end {equation}
for all positive test functions $\gz$ in the space $C_{0}^{2,1}( Q_T)$.\\

If a solution of (\ref{Eq1-1}) exists, it is unique, and we shall denote it by $u_\gm$.
It is not true that problem (\ref{Eq1-1}) can be solved for any positive bounded measure $\gm$ although it is the case if $\gm$ is absolutely continuous with respect to the N-dimensional Hausdorff measure $H^N$.
\bdef {gm}A measure for which the problem can be solved is called a  good measure relative to $g$. The subset of $\frak M_+^1(\Gw)$ of good measures relative to $g$ is denoted by $\CG^{\Gw}(g)$. If  $\gm\in \frak M_+^1(\Gw)$ belongs to $\CG^\Gw(g)$ for any $g$ satisfying (\ref{g1}), is called a universally good measure.
\es
There are many sufficient conditions which insure the solvability of (\ref{Eq1-1}), for example

\begin {equation}\label {Eq2-2}
\dint_{Q_T}g(\BBE[\gm])\gr(x)dx\,dt<\infty,
\end {equation}
where $\BBE[\gm]$ is the heat potential of $\gm$ in $\Gw$, that is the solution $v$ of
\begin {equation}\label {Eq2-3}\left\{\BA {l}
\prt_tv-\Gd v=0\quad\mbox {in }\,Q_T\\[2mm]
\phantom{\prt_tv-\Gd}
v=0 \quad\mbox {in }\,\prt_\ell Q_T\\[2mm]
\phantom{-\Gd}
v(.,0)=\gm\quad\mbox {in }\,\Gw.
\EA \right.
\end {equation}

We recall the parabolic Kato inequality

\blemma{kato} Let $W$ be a domain in $\Gw\ti\BBR$,  $v\in L^1_{loc}(W)$ and
$h\in L^1_{loc}(W)$ such that
\begin{equation}\label{k1}
-\prt_tv+\Gd v\geq h\quad\text {in }\CD'(W).
\end{equation}
Then
\begin{equation}\label{k2}
-\prt_tv_++\Gd v_+\geq h\chi_{[v\geq 0]}\quad\text {in }\CD'(W).
\end{equation}
\es
\Proof Let $\{\gs_j\}$ be a regularizing sequence with compact support in the $N+1 $ ball $\tilde B_{\ge_j}(0)$ ($\ge_j \to 0$ as $j\to\infty$), and
$v_{j}=v\ast \gs_j$. If  $V\subset W$ is such that $\dist (V,W^c)>0$, $v_{j}$ is defined in $V$ whenever $\ge_j <\dist (V,W^c)$. Then
\begin{equation}\label{k3}
-\prt_tv_{j}+\Gd v_{j}\geq h_{j}=h\ast\gs_{j}\quad\text {in }\CD'(V),
\end{equation}
and everywhere in $V$. For $\gd>0$ let
$$j_\gd(r)=\left\{\BA {l}
0\qquad\qquad\text{if }r< -\gd\\[2mm]
\myfrac{(r+\gd)^2}{2\gd}\quad\text{if }-\gd\leq r\leq 0\\[2mm]
r+\myfrac{\gd}{2}\quad\quad\;\text{if }  r> 0. \EA\right.$$
Since
$$-\prt_tj_\gd(v_{j})+\Gd j_\gd(v_{j})
=j'_\gd(v_{j})\left(-\prt_tv_{j}+\Gd v_{j}\right)+j''_\gd(v_{j})
|\nabla v_j|^2\geq j'_\gd(v_{j}) h_{j},
$$
and $\gf\in C^\infty_0(W)$ is nonnegative and has compact support in $V$, it follows that
$$\myint{W}{}j_\gd(v_{j})\left(\prt_t\gf+\Gd\gf\right)dx\,dt
\geq \myint{W}{} j'_\gd(v_{j}) h_{j}\gf\, dx\,dt.
$$
Letting $j\to\infty$, and using the fact that $j_\gd$ and $j'_\gd$ are continuous and, for some subsequence still denoted $\{\ge_j\}$,
$\{(v_{\ge_j},h_{\ge_j})\}$ converges to $(v,h)$ in $L^1_{loc}$ and almost everywhere in $W$, we derive from the Lebesgue theorem
$$\myint{W}{}j_\gd(v)\left(\prt_t\gf+\Gd\gf\right)dx\,dt
\geq \myint{W}{} j'_\gd(v) h\gf\, dx\,dt.
$$
Now $j_\gd (v)$ converges to $v^+$ in $L^1_{loc}$ and $j'_\gd (v(x,t))$ converges to $0$ if $v(x,t)<0$ and to $1$ if $v(x,t)\geq 0$, i.e. to $\chi_{[v\geq 0]}$. Using again the Lebesgue theorem, we obtain
$$\myint{W}{}v^+\left(\prt_t\gf+\Gd\gf\right)dx\,dt
\geq \myint{W}{} \chi_{[v\geq 0]} h\gf\, dx\,dt,
$$
which is (\ref{k2}).\qeda\\

\noindent\Remark In an equivalent way, we can state \rlemma{kato} as follows: {\it If $v\in L^1_{loc}(W)$ and
$h\in L^1_{loc}(W)$ are such that
\begin{equation}\label{k1*}
\prt_tv-\Gd v\leq h\quad\text {in }\CD'(W).
\end{equation}
Then
\begin{equation}\label{k2*}
\prt_tv_+-\Gd v_+\leq h\chi_{[v\geq 0]}\quad\text {in }\CD'(W).
\end{equation}
}

\bdef {tracedef}Let $u\in L_{loc}^1(Q_{T})$. 1- We say that $u$ admits the Radon measure $\gm$ as an initial trace if  it exists
\begin{equation}\label {tr}
\rm ess\,lim_{t\to 0}\myint{\Gw}{}u(.,t)\gf\,dx=\myint{\Gw}{}\gf\,d\gm
\forevery\gf\in C_0(\Gw).
\end {equation}
We shall denote $\gm=Tr_{\Gw}(u)$.\smallskip

\noindent 2-  We say that $u$ admits the outer regular positive Borel measure $\gn\approx (\CS,\gm)$ as an initial trace if  it exists an open subset $\CR\subset\Gw$ and $\gm\in\mathfrak M_+(\CR)$ such that
\begin{equation}\label {trb1}
\rm ess\,lim_{t\to 0}\myint{\Gw}{}u(.,t)\gf\,dx=\myint{\Gw}{}\gf\,d\gm
\forevery\gf\in C_0(\CR).
\end {equation}
and, with $\CS=\Gw\setminus\CR$,
\begin{equation}\label {trb2}
\rm ess\,lim_{t\to 0}\myint{\Gw}{}u(.,t)\gf\,dx=\infty
\forevery\gf\in C_0(\Gw), \gf\geq 0, \gf>0\text { somewhere on }\CS.
\end {equation}
We shall denote $\gn=tr_{\Gw}(u)$.
\es
The trace operator is order preserving. The proof of the following result is straightforward.
\bprop{ord}  Let $u$ and $\tilde u$ in $L^1_{loc}(Q_T)$. \smallskip

\noindent 1- Suppose $Tr_\Gw(u)=\gm$ and $Tr_\Gw(\tilde u)=\tilde\gm$. Then
\begin{equation}\label {ord1}
\tilde u\leq u\Longrightarrow \tilde\gm\leq\gm.
\end {equation}
\noindent 2- Suppose $tr_\Gw(u)=\gn\approx (\CS,\gm)$ and $tr_\Gw(\tilde u)=\tilde\gn\approx(\tilde\CS,\tilde\gm)$. Then
\begin{equation}\label {ord2}
\tilde u\leq u\Longrightarrow \tilde\CS\subset\CS\quad\text{and }\tilde\gm\vline_{\CS^c}
\leq \gm\vline_{\CS^c}.
\end {equation}
\es
The next  classical results characterize the nonnegative supersolutions or subsolutions. We give their proof for the sake of completeness.
\bprop{traceprop1} Let $u\in L^1(Q_{T})$ be a nonnegative supersolution of (\ref{equ}) in $Q_{T}$ such that $g(u)\in L^1(Q_{T})$. Then there exists a positive Radon measure $\gm$ such that $\gm=Tr_{\Gw}(u)$.
\es
\Proof If $0<\gs<\gt<T$ are two Lebesgue points of
$t\mapsto\norm{u(.,t)}_{L^1}$ and $\gf\in C^{2}_0(\Gw)$, $\gf\geq
0$, we set $Q_{\gs,\gt}=\Gw\ti(\gs,\gt)$, take
$\gz(x,t)=\chi_{[\gs,\gt]}(t)\gf(x)$ (by approximations) and
derive from the definition that
\begin {equation}\label {Eq2-1-1}
\myint{\Gw}{}u(.,\gt)\gf \,dx-\myint{\Gw}{}u(.,\gs)\gf \,dx+\dint_{Q_{\gs,\gt}}\left(-u\Gd\gz+\gz\,g(u) \right)dx\,dt\geq0.
\end {equation}
Set
$$H(\gs)=\dint_{Q_{\gs,\gt}}\left(-u\Gd\gf+\gf\,g(u) \right)dx\,dt
$$
Then  $H\in L^1(0,\gt)$ and the mapping
$$\gs\mapsto \Psi(\gs)=\myint{\Gw}{}u(.,\gs)\gf\,dx-H(\gs)
$$
is a.e. nondecreasing on $(0,\gt]$ and it admits an essential limit $L (\gf)\in\BBR$ as $\gs\to 0$. Therefore it exists
$$\ell(\gf)=\rm ess\,lim_{\gs\to 0}\myint{\Gw}{}u(.,\gs)\gf\,dx,
$$
and the mapping $\gf\mapsto \ell (\gf)$ defines a positive Radon measure $\gm$ in $\Gw$.\qeda
\\

 It is possible to get rid of the integrability assumption on $u$ if it is assumed that $u$ vanishes on the boundary and $\Gw$ is bounded.
\bprop{traceprop2} Let $u$ be a positive supersolution of
(\ref{equ}) in $Q_{T}$ which vanishes on $\prt_\ell Q_T$ in the sense that (\ref{Eq2-1''}) holds for all nonnegative $\gz\in C_{\ell,0}^{2,1}(\overline Q_T)$. If $g(u)\in L_{\gr}^1(Q_T)$, there exists $\gm\in\frak M^1_+(\Gw)$ such that $\gm=Tr_{\Gw}(u)$.
\es
\Proof As a test function we take
$\gz(x,t)=\chi_{[\gs,\gt]}(t)\gf_1(x)$ where $\gf_1$ is the first
eigenfunction of $-\Gd$ in $W^{1,2}_0(\Gw)$, $\gf_1\geq 0$ and
$\gl_1$ the corresponding eigenvalue. Thus (\ref{Eq2-1-1}) is
replaced by
\begin {equation}\label {Eq2-1-1'}
\myint{\Gw}{}u(.,\gt)\gf_1 \,dx-\myint{\Gw}{}u(.,\gs)\gf_1 \,dx+\dint_{Q_{\gs,\gt}}\left(\gl_1u+g(u) \right)\gf_1dx\,dt\geq 0.
\end {equation}
If we set
$$X(\gt)=\dint_{Q_{\gs,\gt}}u\gf_1 \,dx\,dt,
$$
and
$$G(\gs)=\dint_{Q_{\gs,\gt}}\gf_1\,g(u)\,dx\,dt,
$$
then (\ref{Eq2-1-1'}) reads as
$$X'(\gs)+\gl_1X(\gs)+G(\gs)\geq X'(\gt)\quad \text{a.e. } 0<\gs<\gt,
$$
which yields to
$$\myfrac{d}{d\gs}\left(e^{\gl_1\gs}X(\gs)-\myint{\gs}{\gt}e^{\gl_1t}\left(G(t)-X'(\gt)\right)dt
\right)\geq 0.$$
The conclusion follows as in \rprop{traceprop1}. Notice also that another choice of test function yields to $u\in L^1(\Gw)$.\qeda\\

For subsolutions of (\ref{equ}) we prove the following.
\bprop{traceprop3} Let $u\in L^1(Q_{T})$ be a nonnegative subsolution of (\ref{equ}) in $Q_{T}$ such that $g(u)\in L^1(Q_{T})$. Then there exists a positive outer regular Borel measure $\gn$ on $\Gw$ such that $\gn=tr_{\Gw}(u)$.
\es
\Proof Defining $H$ as in the proof of \rprop{traceprop1} we obtain that
$$\gs\mapsto \Psi(\gs)=\myint{\Gw}{}u(.,\gs)\gf\,dx+H(\gs)
$$
is nonincreasing on $(0,\gt]$ and it admits a limit $L^* (\gf)\in (-\infty,\infty]$ as $\gs\to 0$. For any $\xi\in\Gw$ the following dichotomy holds,\smallskip

\noindent (i) either there exists a $\gf\in C^{2}_0(\Gw)$ verifying $\gf(\xi)>0$ such that $L(\gf)<\infty$,\smallskip

\noindent (ii) or for any $\gf\in C^{2}_0(\Gw)$ verifying $\gf(\xi)>0$, $L(\gf)=\infty$.\smallskip

\noindent The set $\CR(u)$ of $\xi$ such that (i) occurs is open and
there exists $\gm\in \mathfrak M_+(\CR(u))$ such that
$$L(\gf)=\myint{\CR(u)}{}\gf\,d\gm \forevery \gf\in C_0(\CR(u)).
$$
The set $\CS(u)=\Gw\setminus\CS(u)$ is relatively closed in $\Gw$. Further, if $\gf\in C_0(\Gw)$ is nonnegative and positive somewhere on $\CS(u)$, there holds
$$\rm ess\,lim_{t\to 0}\myint{\Gw}{}\,u(.,\gs)\gf\,dx=\infty.
$$
The outer regular Borel measure $\gn$ is defined for any Borel subset
$E\subset\Gw$ by
$$\gn(E)=\left\{\BA {l}\myint{E}{}\,d\gm\quad \text {if }E\subset \CR(u)\\
\infty \quad\; \;\quad \text {if }E\cap\CS(u)\neq\emptyset .\EA\right.
$$
\qeda\\

The next lemma is the parabolic counterpart of an elliptic result proved in \cite{BP}
\blemma {n-lim}Let $f\in L_\gr^1(Q_T)$ and $u\in L^1(Q_T)$ such that
\begin {equation}\label{lin1}
\dint_{Q_T}u(\prt_t\gz+\Gd\gz)dx\,dt=-\dint_{Q_T}f\gz\,dx\,dt
\end {equation}
for every $\gz\in C^{2,1}_{\ell,0}(\overline Q_T)$. Then
\begin {equation}\label {lin2}
\lim_{n\to\infty}n\dint_{Q_T\cap\{\gr(x)\leq n^{-1}\}}\abs udx\,dt=0.
\end {equation}
\es
\Proof We assume first that $f\geq 0$, then $u\geq 0$. Let $H$ be a nondecreasing concave $C^2$ function such that $H(0)=0$, $H''(t)=-1$ for $0\leq t\leq 1$ and $H(t)=1$ for $t\geq 2$. Let $\xi_0$ be the solution of
\begin{equation}\label{xi-0}\left\{\BA {l}
\prt_t\xi_0+\Gd\xi_0=-1\quad\text{in }Q_T\\[2mm]
\phantom{\prt_t\xi_0}
\xi_0(.,T)=0\quad\text{in }\bar\Gw\\[2mm]
\phantom{\prt_t\xi_0}
\xi_0(x,t)=0\quad\text{in }\prt\Gw\ti[0,T].
\EA\right.\end{equation}
Let $w_n=n^{-1}H(n\xi_0)$, then
$$-\prt w_n-\Gd w_n\geq -nH''(n\xi_0)\abs{\nabla \xi_0}^2\geq n\chi_{_{\{\xi_0\leq n^{-1}\}}}\abs{\nabla \xi_0}^2.
$$
Therefore
$$\dint_{Q_T}fw_ndx\,dt=-\dint_{Q_T}u(\prt_t w_n+\Gd w_n)dx\,dt\geq n
\dint_{Q_T}\abs{\nabla \xi_0}^2udx\,dt.
$$
But $w_n\leq\min\{\xi_0,n^{-1}\}$, therefore, by the Lebesgue theorem,
$$0=\lim_{n\to\infty}\dint_{Q_T}fw_ndx\,dt=
\lim_{n\to\infty}n\dint_{Q_T}\abs{\nabla \xi_0}^2udx\,dt.
$$
Let $\ge>0$, by Hopf lemma on $Q_{T-\ge}$, there exists $c_1>0$, $c_2>0$ such that $\abs{\nabla \xi_0}\geq c_1$ on $\prt\Gw\ti[0,T-\ge]$; thus $c_2\xi_0\leq\gr\leq c_2^{-1}\xi_0$ and
$$\lim_{n\to\infty}n\dint_{Q_{T-\ge}\cap\{\xi_0(x)\leq n^{-1}\}}udx\,dt=0.
$$
Clearly we can extend $f$ to be zero for $t>T$ and $\tilde u$ to be the weak solution of
$$\left\{\BA {l}
\prt_t\tilde u+\Gd \tilde u=0\quad\text{in }Q_{T,T+\ge}\\[2mm]
\phantom{\prt_t\xi_0}
\tilde u(.,T)=u(.,T)\quad\text{in }\bar\Gw\\[2mm]
\phantom{\prt_t\xi_0}
\tilde u(x,t)=0\quad\text{in }\prt\Gw\ti[T,T+\ge].
\EA\right.$$
Notice that it is always possible to assume that $T$ is a Lebesque point of
$t\mapsto\norm{u(.,t)}_{L^1}$ inasmuch this function is actually continuous. Replacing $T$ by $T+\ge$, we derive (\ref{lin2}).
Next, if $u$ has not constant sign, we denote by $v$ the weak solution of
$$\left\{\BA {l}
\prt_t v-\Gd v=\abs f \quad\text{in }Q_T\\[2mm]
\phantom{\prt_t\xi_0}
v(.,0)=0\quad\text{in }\bar\Gw\\[2mm]
\phantom{\prt_t\xi_0}
v(x,t)=0\quad\text{in }\prt\Gw\ti[0,T].
\EA\right.$$
Then $\abs u\leq v$ and the proof follows from the first case.
\qeda\\
\blemma {n-lim2}Let $f\in L_\gr^1(Q_T)$ and $u\in L^1(Q_T)$ such that
\begin {equation}\label{lin3}
-\dint_{Q_T}u(\prt_t\gz+\Gd\gz)dx\,dt\leq\dint_{Q_T}f\gz\,dx\,dt
\end {equation}
for every $\gz\in C^{2,1}_{\ell,0}(\overline Q_T)$, $\gz\geq 0$. Then, for the same class of test functions $\gz$, there holds
\begin {equation}\label {lin4}
-\dint_{Q_T}(\prt_t\gz+\Gd\gz)u_+dx\,dt\leq\dint_{Q_T\cap\{u\geq 0\}}f\gz\,dx\,dt.
\end {equation}
\es
\Proof By \rlemma{kato}, (\ref{lin4}) holds for any $\gz\in
C^{2,1}_{0}( Q_T)$. Let $\{\gg_n\}$ be a sequence of functions in
$C^{2,1}_{0}( Q_T)$ such that $0\leq\gg_n\leq1$, $\gg_n(x,t)=1$ if
$\gr(x)\geq n^{-1}$ or $t\geq n^{-1}$,
$\norm{\nabla{\gg_n}}_{L^\infty}\leq Cn$,
$\norm{\Gd{\gg_n}}_{L^\infty}\leq Cn^2$ and
$\norm{\prt_t{\gg_n}}_{L^\infty}\leq Cn$. Given $\gz\geq 0$ in
$C^{2,1}_{\ell,0}(\overline Q_T)$, $\gz\gg_n$ is an admissible
test function for Kato's inequality (\ref{lin4}), thus
\begin {equation}\label {lin5}
-\dint_{Q_T}(\prt_t(\gz\gg_n)+\Gd(\gz\gg_n))u_+dx\,dt\leq\dint_{Q_T\cap\{u\geq 0\}}f(\gz\gg_n)\,dx\,dt.
\end {equation}
When $n\to\infty$ the right-hand side of (\ref{lin5}) converges to the right-hand side of (\ref{lin4}). Moreover
$\prt_t(\gz\gg_n)=\gg_n\prt_t\gz+\gz\prt_t\gg_n$, $\nabla(\gz\gg_n)=\gg_n\nabla\gz+\gz\nabla\gg_n$  and $\Gd(\gz\gg_n)=\gg_n\Gd\gz+\gz\Gd\gg_n+2\nabla\gz.\nabla\gg_n$. Thus
$$\prt_t(\gz\gg_n)+\Gd(\gz\gg_n)=\gg_n\prt_t\gz+\gz\prt_t\gg_n+
\gg_n\Gd\gz+\gz\Gd\gg_n+2\nabla\gz.\nabla\gg_n.
$$
Since $\gz$ vanishes on $\prt\Gw\ti[0,T]$ and is bounded with bounded gradient, there holds
$$\abs{\dint_{Q_T}\left(\gz\prt_t\gg_n+\gz\Gd\gg_n+2\nabla\gz.\nabla\gg_n\right)u^+\,dx\,dt}\leq Cn\dint_{Q_T\cap\{\gr(x)\leq n^{-1}\}}u^+
\,dx\,dt$$
which goes to 0 as $n\to\infty$. This implies (\ref{lin4}).
\qeda\\

If we deal with subsolution or supersolutions of problem (\ref{Eq1-1}) we have the following results
\bth{trsub}
Let $\gm\in\frak M_+^1(\Gw)$ and $u$ be a nonnegative subsolution of (\ref{Eq1-1}). Then
 the initial trace of $u$ is a positive Radon measure $\tilde \gm$ such that
$\tilde \gm\leq \gm$. Furthermore, if (\ref{Eq1-1}) admits a weak solution $u_\gm$, there holds
 $u\leq u_\gm$.
\es
\Proof {\it Step 1. There holds $\tilde \gm\leq \gm$.}  If $\gs$ is a Lebesgue point of $t\mapsto\norm{\tilde u(.,t)}_{L^1}$ and $\gf\in C^{2}_0(\Gw)$, $\gf\geq 0$, we can take $\gz(x,t)=\chi_{[0,\gs]}(t)\gf(x)$ (by approximations) and derive from (\ref{Eq2-1'}) that
\begin {equation}\label {Eq2*-1-1}
\myint{\Gw}{} u(.,\gs)\gf \,dx+\dint_{Q_{\gs}}\left(-u\Gd\gz+\gz\,g( u) \right)dx\,dt\leq\myint{\Gw}{}\gf\,d\gm,
\end {equation}
thus, by \rprop {traceprop3}, using the fact that $u\in L^1(Q_T)$ and
$g(u)\in L^1_\gr(Q_T)$,
\begin {equation}\label {Eq2'-7}\rm ess\,\lim_{\gs\to 0}
\myint{\Gw}{} u(.,\gs)\gf \,dx\leq \myint{\Gw}{}\gf\,d\gm.
\end {equation}
It follows that the initial trace $\tilde\gn\approx(\CS(u),\tilde\gm$ has no singular part ($\CS(u)=\emptyset$) and $\tilde \gm\leq \gm$.
This implies that $\gf\mapsto m(\gf)$ is a measure dominated by $\gm$ that we shall denote by $\tilde \gm$. It represents the initial trace of $\tilde u$, and we shall denote it by
\begin {equation}\label {Eq2''-7}
\tilde \gm=Tr_\Gw(\tilde u).
\end {equation}
Next we take $\gz\in C^{2,1}_{\ell,0}(\overline Q_T)$, $\gz\geq 0$, and get at any Lebesgue point $\gs$ as in \rprop{traceprop1}-\rprop{traceprop3}
$$
\dint_{Q_{\gs,T}}\left(-u\prt_t\gz-u\Gd\gz+\gz\,g( u) \right)dx\,dt\leq
\myint{\Gw}{} u(.,\gs)\gz \,dx,
$$
we derive, by letting $\gs\to 0$,
\begin {equation}\label {L}
\dint_{Q_{T}}\left(-u\prt_t\gz-u\Gd\gz+\gz\,g( u) \right)dx\,dt\leq
\myint{\Gw}{} \gz(.,0)\,d\tilde\gm.
\end {equation}
{\it Step 2. There exists $u_{\tilde\gm}$ and $u_{\tilde\gm}\leq u_\gm$.} For $k>0$ set $g_k(r)=\min\{g(r),k\}$ and let $u=u^k_{\tilde\gm}$ be the solution of
\begin {equation}\label {Eq1-epsilon}\left\{\BA {l}
\prt_tu-\Gd u+g_k(u)=0\quad \text {in } Q_T:=\Gw\ti(0,T)\\[2mm]
\phantom{\prt_t-\Gd u+g_k(u)}
u=0\quad \text {in } \prt_\ell Q_T:=\prt\Gw\ti(0,T)\\[2mm]
\phantom{\Gd u+g_k(u)}
u(.,0)=\tilde\gm\quad \text {in } \Gw.
\EA\right.
\end {equation}
Defining in the same way $u^k_{\gm}$, we obtain $u^k_{\tilde\gm}\leq u^k_{\gm}$, $u^k_{\tilde\gm}\geq u^{k'}_{\tilde\gm}$ and
$u^k_{\gm}\geq u^{k'}_{\gm}\geq u_\gm$ for $k'>k>0$.
If $\gz\in C_{\ell,0}^{2,1}(\overline Q_T)$ is nonnegative, there holds
\begin{equation}\label{cv}
\dint_{Q_T}\gz\,g_k(u^k_\gm) dx\,dt=\myint{\Gw}{}\gz\,d\gm+
\dint_{Q_T}\left(\prt_t\gz+\Gd\gz \right)u^k_\gm dx\,dt.
\end{equation}
Clearly $u^k_\gm$ converges to some $U\geq u_\gm$ when $k\to\infty$, the right-hand side of (\ref{cv}) converges to
$$\myint{\Gw}{}\gz\,d\gm+
\dint_{Q_T}\left(\prt_t\gz+\Gd\gz \right)U dx\,dt,
$$
and $g_k(u^k_\gm)$ converges to $g(U)$ a. e. By Fatou
$$\dint_{Q_T}\gz\,g(U) dx\,dt\leq \liminf_{k\to\infty}\dint_{Q_T}\gz\,g_k(u^k_\gm) dx\,dt,
$$
thus, using the monotonicity of $g$,
\begin{equation}\label{cv2}
\dint_{Q_T}\gz\,g(u_\gm) dx\,dt\leq \dint_{Q_T}\gz\,g(U) dx\,dt\leq \myint{\Gw}{}\gz\,d\gm+
\dint_{Q_T}\left(\prt_t\gz+\Gd\gz \right)U dx\,dt.
\end{equation}
Because $u_\gm$ satisfies (\ref{Eq2-1}), all the three terms in
(\ref{cv2}) are equal, $U=u_\gm$ and
\begin{equation}\label{cv3}
\lim_{k\to\infty}\dint_{Q_T}\gz\,g_k(u^k_\gm) dx\,dt=\dint_{Q_T}\gz\,g(u_\gm) dx\,dt.
\end{equation}
Next $u^k_{\tilde\gm}$ decreases and converges to some $\tilde U$,
 $g_k(u^k_{\tilde\gm})\to g(\tilde U)$ a.e., and
\begin{equation}\label{cv4}
\dint_{Q_T}\gz\,g(\tilde U) dx\,dt\leq \lim_{k\to\infty}\dint_{Q_T}\gz\,g_k(u^k_{\tilde\gm}) dx\,dt=
\myint{\Gw}{}\gz\,d\tilde\gm+
\dint_{Q_T}\left(\prt_t\gz+\Gd\gz \right)\tilde U dx\,dt.
\end{equation}
Since $0\leq \gz\,g_k(u^k_{\tilde\gm})\leq \gz\,g_k(u^k_\gm)$. In order to prove that
\begin{equation}\label{cv5}
 \lim_{k\to\infty}\dint_{Q_T}\gz\,g_k(u^k_{\tilde\gm}) dx\,dt
=\dint_{Q_T}\gz\,g(\tilde U) dx\,dt,
\end{equation}
we use the following classical result : {\it Let $h_n\geq\tilde h_n\geq 0$ two sequences of measurable functions in some measured space $(G,\Gs,dm)$ which converge a. e. in $G$ to $h$ and $\tilde h$ respectively. Then}
$$\lim_{n\to\infty}\myint{G}{}h_ndm=\myint{G}{}hdm\Longrightarrow
\lim_{n\to\infty}\myint{G}{}\tilde h_ndm=\myint{G}{}\tilde hdm.
$$
Therefore (\ref{cv3}) implies (\ref{cv5}). From (\ref{cv4}) we get
\begin{equation}\label{cv6}
\dint_{Q_T}\gz\,g(\tilde U) dx\,dt=
\myint{\Gw}{}\gz\,d\tilde\gm+
\dint_{Q_T}\left(\prt_t\gz+\Gd\gz \right)\tilde U dx\,dt.
\end{equation}
This relation is valid with any $\gz\in C_{\ell,0}^{2,1}(\overline Q_T)$ with constant sign. It implies in particular that $Tr_\Gw(\tilde U)=\tilde \gm$. Thus $u_{\tilde\gm}$ exists and $\tilde U=u_{\tilde\gm}$.\smallskip

\noindent{\it Step 3. We claim that $u\leq u_{\tilde\gm}$.} Set
$w=u-u_{\tilde\gm}$, it follows from (\ref{L}),
\begin{equation}\label{cv7}
\dint_{Q_T}\left(-w\prt_t\gz-w\Gd\gz +(g(u)-g(u_{\tilde\gm}))\gz\right) dx\,dt
\leq 0
\end{equation}
for any $\gz\in C^{2,1}_{\ell,0}(\overline Q_T)$, $\gz\geq 0$. Using \rlemma {n-lim2} we derive
\begin{equation}\label{cv8}
\dint_{Q_T\cap\{w^+\geq 0\}}\left(-(\prt_t\gz-\Gd\gz)w_++ +(g(u)-g(u_{\tilde\gm}))\gz\right) dx\,dt
\leq 0\end{equation}
We take $\gz=\xi_0$ given by (\ref{xi-0}). Since $g$ is nondecreasing, we derive
\begin{equation}\label{cv9}
\dint_{Q_T\cap\{w^+\geq 0\}}w_+ dx\,dt
\leq 0.\end{equation}
Thus $u\leq u_{\tilde\gm}\leq u_\gm$.
\qeda
\\

\noindent\Remark It is noticeable that Step-2 of the proof of \rth{trsub} can be stated in the following way. {\it If $\gm\in\mathfrak M_+^1(\Gw)$ is a good measure, any measure $\tilde\gm$ such that $0\leq\tilde\gm\leq\gm$ is a good measure.}\\

Consider $\gm\in\mathfrak M^1_+(\Gw)$. The relaxation phenomenon associated to (\ref{Eq1-1}) can be constructed in the following way. Let $\{g_k\}$ be an increasing sequence of continuous nondecreasing functions defined on $\BBR$, vanishing on $(-\infty,0]$ and such that
\begin {equation}\label {Eq2-4}\BA {l}
(i)\qquad 0\leq g_k(r)\leq c_kr^{p}+c'_k\forevery r\geq 0,\;\;\forall k>0\\[2mm]
(ii)\qquad\displaystyle\lim_{k\to\infty} g_k(r)= g(r) \forevery r\in\BBR,
\EA\end {equation}
for some positive constants $c_k$ and $c'_k$ and $p\in(1,(N+2)/(N+1))$. Since (\ref{Eq2-2}) is satisfied, there exists a unique solution $u=u_k$ to
\begin {equation}\label {Eq2-5}\left\{\BA {l}
\prt_tu-\Gd u+g_k(u)=0\quad \text {in } Q_T\\[2mm]
\phantom{\prt_t-\Gd u+g_k(u)}
u=0\quad \text {in } \prt_\ell Q_T\\[2mm]
\phantom{\Gd u+g_k(u)}
u(.,0)=\gm\quad \text {in } \Gw.
\EA\right.
\end {equation}

It is noticeable that, if the assumption $\gm\in\mathfrak M^1_+(\Gw)$ were replaced by $\gm\in\mathfrak M^0_+(\Gw)$, the exponent $p$ in (\ref{Eq2-4}) should have been taken smaller than $(N+2)/N$. In the sequel $C$ will denote a positive constant, depending on the data, not on $k$, the value of which may change from one occurrence to another. Our first result points out the relaxation phenomenon associated to the sequence $\{u_k\}$.
\bth {redmes} When $k\to\infty$, the sequence $\{u_k\}$ converges in $L^1(Q_T)$ to a some nonnegative function $u^*$ such that
$g(u^*)\in L_\gr^1(Q_T)$, and there exists a positive measure $\gm^*$ smaller that $\gm$ with the property that
\begin {equation}\label {Eq2-5*}\left\{\BA {l}
\prt_tu^*-\Gd u^*+g(u^*)=0\quad \text {in } Q_T\\[2mm]
\phantom{\prt_t-\Gd u^*+g(u^*)}
u^*=0\quad \text {in } \prt_\ell Q_T\\[2mm]
\phantom{\Gd u^*+g(u^*)}
u^*(.,0)=\gm^*\quad \text {in } \Gw.
\EA\right.
\end {equation}
Furthermore $u^*$ is the largest subsolution of problem (\ref{Eq1-1}).
\es
\Proof By \cite [Lemma1.6] {MV1} there holds
\begin {equation}\label {Eq2-6}
\norm {u_k}_{L^1}+\norm {g_k(u_k)}_{L_\gr^1}\leq C\myint{\Gw}{}\gr\,d\gm,
\end {equation}
and, by the maximum principle,
\begin {equation}\label {Eq2-7}
u_k\leq \BBE[\gm]\quad\text{in }\;Q_T.
\end {equation}
For any $\ge>0$ we denote $Q_{\ge,T}=\Gw\ti[\ge,T]$. Since $\BBE[\gm]$ is uniformly bounded in $Q_{\ge,T}$ for any $\ge>0$, it follows by the parabolic equations regularity theory that,  $u_k$ is bounded in $C^{1+\ga,\ga/2}(Q_{\ge,T})$ for any $0<\ga<1$. Furthermore, if $k'>k$, $g_{k'}(u_k)\geq g_{k}(u_k)$ thus $u_k$ is a super-solution for the equation satisfied by
$u_{k'}$. This implies $u_k\geq u_{k'}$ and $u^*:=\lim_{k\to\infty}u_k$ exists and satisfies
$$u^*\leq \BBE[\gm]\quad\text{in }\;Q_T.
$$
Because of (\ref {Eq2-7}) uniform boundedness holds also in $L^p(Q_T)$, for any $p\in [1,(N+2)/(N+1))$. By the Lebesgue theorem the convergence occurs in $L^p(Q_T)$ too, for any $p\in [1,(N+2)/(N+1))$, and locally uniformly in $Q_T$ by the standard regularity theory. By continuity $g_k(u_k)$ converges to $g(u^*)$ uniformly in $Q_{\ge,T}$, thus $u^*$ satisfies
$$\prt_tu^*-\Gd u^*+g(u^*)=0\quad \text {in } Q_T
$$
and vanishes on $\prt_\ell Q_T$. By the Fatou theorem
$$\dint_{Q_T} g(u^*)\gz \,dx\,dt\leq\liminf_{k\to\infty}\dint_{Q_T} g_k(u_k)\gz \,dx\,dt,
$$
for any $\gz\in C(\overline Q_T)$, $\gz\geq 0$, and there exists a positive measure
$\gl$ in $Q_T$ such that
$$g_k(u_k)\rightarrow g(u^*)+\gl,
$$
weakly in the sense of measures. Thus for any
 $\gz\in C_{\ell,0}^{2,1}(\overline Q_T)$, there holds
\begin{equation}\label{Eq2-8}
\dint_{Q_T}\left(-u^*\prt_t\gz-u^*\Delta \gz+g(u^*)\gz\right)dx\,dt
=\myint{\Gw}{}\gz(x,0)\, d\gm-\dint_{Q_T}\gz\, d\gl.
\end{equation}
Since $g_k(u_k)$ converges to $g(u^*)$ uniformly in $Q_{\ge,T}$ for any $\ge>0$, the measure $\gl$ is concentrated on $\overline\Gw\ti\{0\}$. We denote by $\tilde\gl$ its restriction to
$\Gw\ti\{0\}$, set
$$\gm^*=\gm-\tilde\gl,
$$
and derive from (\ref{Eq2-8}),
\begin{equation}\label{Eq2-9}
\dint_{Q_T}\left(-u^*\prt_t\gz-u^*\Delta \gz+g(u^*)\gz\right)dx\,dt
=\myint{\Gw}{}\gz(x,0)\, d\gm^*.
\end{equation}
This implies $u^*=u_{\gm^*}$ and $Tr_\Gw(u^*)=\gm^*$,
thus $\gm^*$ is a positive measure. Let $v$ be a nonnegative subsolution of problem (\ref{Eq2-1}). By \rprop{traceprop3} there exists $\tilde\gm\in \mathfrak M^1_+(\Gw)$ such that
$Tr_\Gw(v)=\tilde\gm$ and $\tilde\gm\leq\gm$. Since $g_k(v)\leq g(v)$, $u$ is a subsolution for problem (\ref{Eq2-5}). By \rth{trsub} $v\leq u_k:= u_{k,\gm}$. Thus $\lim_{k\to\infty}u_k=u^*\geq v$.\qeda

\bth {redmes2} The reduced measure $\gm^*$ is the largest good
measure smaller than $\gm$. \es
\Proof Clearly $\gm^*$ is a good measure smaller than $\gm$. Assume now that $\tilde\gm$ is a good measure smaller than $\gm$. Then $u_{\tilde\gm}$ is a subsolution for problem (\ref{Eq2-1}). By (\rth{trsub}) $u_{\gm^*}$ is larger than $u_{\tilde\gm}$. Thus
$Tr_\Gw(u_{\tilde\gm})=\tilde\gm\leq Tr_\Gw(u_{\gm^*})=\gm^*$.\qeda \\

The next technical result characterizes the good measures
\bth {redmes3} Let $\gm\in \mathfrak M_+(\Gw)$. Then $\gm\in\CG^\Gw(g)$ if and only if $g_k(u_k)\to g(u)$ in the weak sense of measures in $\frak M ^1(Q_T)$.
\es
\Proof Assume $g_k(u_k)\to g(u)$ in the weak sense of measures in $\frak M ^1(Q_T)$. Letting $k\to\infty$ in (\ref{cv}), we obtain (\ref{Eq2-1}) for any $\gz\in C^{2,1}_{\ell,0}(\overline Q_T)$. Thus
$u^*=u_\gm$. Thus $\gm^*=\gm$ and $\gm$ is a good measure. Conversely, assume $\gm$ is a good measure. By \rth{redmes2}, $\gm^*=\gm$. Thus $u_k\to u^*=u_{\gm}$ and $u_k\to u_{\gm}$ in $L^1(\Gw)$ and a.e. in $\Gw$. Assume $\gz\in C^{2,1}_{\ell,0}(\overline Q_T)$, $\gz\geq 0$. We let $k\to\infty$ in (\ref{cv}) and derive
\begin{equation}\label{cv3*}
\lim_{k\to\infty}\dint_{Q_T}g_k(u_k)\gz\,dx\,dt=\int_{\Gw}\gz\,d\gm+\dint_{Q_T}\left(\prt\gz+\Gd\gz\right)u_\gm\,dx\,dt= \dint_{Q_T}g(u_\gm)\gz\,dx\,dt,
\end {equation}
by (\ref{Eq2-1}). Because $\{g_{k}(u_{k)}\}$ is uniformly bounded in $L^1_{\gr}(Q_{T})$, the result follows by density.\qed\\

As in \cite{BP} an easy consequence of \rth {redmes2} is the following result which points out the fact that $\gm$ and $\gm^*$ differ only on a set with zero N-dimensional Hausdorff measure.
\bcor{redmes4} Let $\gm\in \frak M^1_+(\Gw)$. There exists a Borel set $E\subset\Gw$, with Hausdorff measure $H^N(E)=0$, such that $(\gm-\gm^*)(E^c)=0$.
\es
\Proof Let $\gm=\gm_r+\gm_s$ be the Lebesgue decomposition of
$\gm$, $\gm_r$ (resp. $\gm_s$) being the absolutely continuous (resp. singular) part relative to the Hausdorff measure $H^N$ in $\BBR^N$. Both measures are positive. Since $\gm_r\in L^1_\gr(\Gw)$, it is a good measure. Then
$\gm_r\leq \gm^*$ by \rth{redmes2}. Therefore
$$0\leq \gm-\gm^*\leq \gm-\gm_r=\gm_s.
$$
Since $\gm_s$ is singular relative to $H^N$, its support $E$ satisfies
 $H^N(E)=0$. This implies the claim.\qeda\\

 \bcor{redmes5} Let $\gm\in \frak M^1_+(\Gw)$ such that $\gm(E)=0$ for any Borel set $E\subset\Gw$ with $H^N(E)=0$. Then $\gm$ is a good measure.
\es
\Proof Let $E\subset\Gw$ is a Borel set with $H^N(E)=0$, then $\gm_r(E)=0$. Since $\gm(E)=0$, it implies $\gm_s(E)=0$. Because the support of $\gm_s$ is a set with zero N-dimensional Hausdorff, $\gm=\gm_r=\gm^*$.\qeda
\bth{ord} Let $\gm_1, \gm_2\in\frak M_+^1(\Gw)$. If $\gm_1\leq\gm_2$, then $\gm_1^*\leq \gm_2^*$. Furthermore
\begin{equation}\label{cont}
 \gm_2^*-\gm_1^*\leq \gm_2-\gm_1.
\end {equation}
\es
\Proof For $k>0$ let $u=u_{k,i}$ ($i=1,2$) be the solution of
\begin {equation}\label {approx}\left\{\BA {l}
\prt_tu-\Gd u+g_k(u)=0\quad \text {in } Q_T\\[2mm]
\phantom{\prt_t-\Gd u+g_k(u)}
u=0\quad \text {in } \prt_\ell Q_T\\[2mm]
\phantom{\Gd u+g_k(u)}
u(.,0)=\gm_i\quad \text {in } \Gw.
\EA\right.
\end {equation}
Since $\gm_1\leq\gm_2$, $u_{k,1}\leq u_{k,2}$. By the convergence result of \rth{redmes}, the relaxed solutions $u_{i}^*$ satisfies $u_{1}^*\leq u_{2}^*$. Since $\gm_i^*=Tr_\Gw(u_{i}^*)$, it follows $\gm_1^*\leq \gm_2^*$. We turn now to the proof of (\ref{cont}). If $\gz\in C^{2,1}(\bar Q_{t})$, $\gz\geq 0$, which vanishes on $\prt_{\ell}Q_{t}$, we have from the weak formulation
$$\BA {l}\dint_{{Q_{t}}}\left(-(u_{k,2}-u_{k,1})(\prt_{t}\gz+\Gd\gz)+\gz(g_{k}(u_{2}^*)-g_{k}(u_{2}^*))\right)dx\,dt\\[2mm]
\phantom{----------}
=\myint{\Gw}{}\gz(x,0) d(\gm_2-\gm_1)-\myint{\Gw}{}\gz(x,t) (u_{k,2}-u_{k,1})dx
\EA$$
We fix $\xi\in C^{2}_{0}(\bar\Gw)$, $\xi\geq 0$ and choose for $\gz$ the solution of
$$\left\{\BA {l}\prt_{t}\gz+\Gd\gz=0\qquad\text{in }Q_{t}\\[2mm]
\phantom{\prt_{t}\gz+\Gd}
\gz=0\qquad\text{on }\prt_{\ell}Q_{t}\\[2mm]
\phantom{\prt_{t}\gz}
\gz(x,t)=\xi\qquad\text{in }\Gw,
\EA\right..$$
Then, letting $k\to\infty$, we derive
$$\myint{\Gw}{}(u_{2}^*-u_{1}^*)(x,t)\xi dx\leq \myint{\Gw}{}\gz(x,0) d(\gm_2-\gm_1).
$$
Finally, if $t\to 0$, using the trace property and the fact that  $\gz(x,0)\to \xi$ in $C_{0}(\bar\Gw)$, we obtain
$$\myint{\Gw}{}\xi d(\gm^*_2-\gm^*_1)\leq \myint{\Gw}{}\xi d(\gm_2-\gm_1).
$$
This implies (\ref{cont}).
\qeda
\bcor {smlrcor}If $\gm$ is a good measure, any positive measure $\gn$ smaller than $\gm$ is a good measure.
\es
\Proof Let $\gn\in \frak M^1_{+}(\Gw)$, $\gn\leq \gm$. By (\ref{cont})
$$0\leq \gn-\gn^*\leq \gm-\gm^*.
$$
Thus $\gm=\gm^*\Longrightarrow\gn=\gn^*$.\qeda
 \bcor{redmes4'} Let $\gm_{1},\gm_{2}\in \frak M^1_+(\Gw)$. 1- If $\gm_{1}$ and $\gm_{2}$ are good measures, then so is $\inf\{\gm_{1}, \gm_{2}\}$ and $\sup\{\gm_{1}, \gm_{2}\}$.\smallskip

 \noindent 2- If $E\subset\Gw$ is a Borel set and $\gm\in\frak M_{+}^1\Gw)$, $\gm^*\vline_{E}=[\gm\vline_{E}]^*$
 \smallskip

 \noindent 3- Assume that $\gm_{1}$ and $\gm_{2}$ are mutually singular. Then $(\gm_{1}+\gm_{2})^*=\gm_{1}^*+\gm_{2}^*$.
\es
  \Proof 1- The fact that $\inf\{\gm_{1}, \gm_{2}\}$ is a good measure is clear from \rcor{smlrcor}. Let $\gn=\sup\{\gm_{1}, \gm_{2}\}$. Then $\gm_{1}\leq \gn^*$ and $\gm_{2}\leq \gn^*$. Then
  $\gn=\sup\{\gm_{1}, \gm_{2}\}\leq\gn^*$.\smallskip

  \noindent 2- We recall that $\gm\vline_{E}(A)=\gm(E\cap A)$, for any Borel subset $A$ of $\Gw$. We can also write $\gm\vline_{E}=\chi_{_{E}}\gm$. Since $\gm\geq \gm^*$, $\chi_{_{E}}\gm\geq \chi_{_{E}}\gm^*$ and also $\gm^*\geq \chi_{_{E}}\gm^*$. Thus $\chi_{_{E}}\gm^*$ is a good measure and
  $[\chi_{_{E}}\gm]^*\geq \chi_{_{E}}\gm^*$ by \rth{redmes2}.
  Conversely, $[\chi_{_{E}}\gm]^*\leq \chi_{_{E}}\gm$ implies that $\chi_{_{E}}[\chi_{_{E}}\gm]^*=[\chi_{_{E}}\gm]^*$. But
  $\chi_{_{E}}\gm\leq \gm$ implies $[\chi_{_{E}}\gm]^*\leq \gm^*$ and therefore
  $[\chi_{_{E}}\gm]^*=\chi_{_{E}}[\chi_{_{E}}\gm]^*\leq \chi_{_{E}}\gm^*$.
  \smallskip

 \noindent 3- If $\gm_{1}$ and $\gm_{2}$ are mutually singular, then so are $\gm^*_{1}$ and $\gm^*_{2}$. Actually, $\gm_{1}+\gm_{2}=\sup\{\gm_{1},\gm_{2}\}$ and $\gm^*_{1}+\gm^*_{2}=\sup\{\gm^*_{1},\gm^*_{2}\}$. By assertion 1,  $[\sup\{\gm^*_{1},\gm^*_{2}\}]^*=\sup\{\gm^*_{1},\gm^*_{2}\}$. Then $\gm^*_{1}+\gm^*_{2}$ is a good measure smaller than $\gm_{1}+\gm_{2}$, thus $\gm^*_{1}+\gm^*_{2}\leq (\gm_{1}+\gm_{2})^*$.
Conversely, there exist two disjoint Borel sets $A$ and $B$ such that
$\gm_{1}=\chi_{_{A}}\gm_{1}$ and $\gm_{2}=\chi_{_{B}}\gm_{2}$ and $\gm_{1}+\gm_{2}=\chi_{_{A}}\gm_{1}+\chi_{_{B}}\gm_{2}$. Thus
 $(\gm_{1}+\gm_{2})^*=(\chi_{_{A}}\gm_{1}+\chi_{_{B}}\gm_{2})^*$ and
 $\chi_{_{A}}(\gm_{1}+\gm_{2})^*=(\chi_{_{A}}\gm_{1}+\chi_{_{A}}\gm_{2})^*=\chi_{_{A}}\gm^*_{1}=\gm^*_{1}$. Similarly, $\chi_{_{B}}(\gm_{1}+\gm_{2})^*=(\chi_{_{B}}\gm_{1}+\chi_{_{B}}\gm_{2})^*=\chi_{_{B}}\gm^*_{2}=\gm^*_{2}$. Since
 $$(\gm_{1}+\gm_{2})^*=\chi_{_{A\cup B}}(\gm_{1}+\gm_{2})^*=\chi_{_{A}}(\gm_{1}+\gm_{2})^*+\chi_{_{B}}(\gm_{1}+\gm_{2})^*,$$
 the result follows.\qeda

\bth{structh}The set $\CG^\Gw(g)$ is a convex lattice. Furthermore
\begin {equation}\label{inf}
[\inf\{\gm,\gn\}]^*=\inf\{\gm^*,\gn^*\},
\end {equation}
and
\begin {equation}\label{sup}
[\sup\{\gm,\gn\}]^*=\sup\{\gm^*,\gn^*\}.
\end {equation}
\es
\Proof For the sake of completeness, we present the proofs of
these assertions which actually the ones already given in
\cite{BMP}. Let $\gm_1,\gm_2\in \CG^\Gw(g)$ and
$\gn=\sup\{\gm_1,\gm_2\}$. Since $\gm_i\leq \gn$, it follows from
\rth{ord} that $\gm_i=\gm^*_i\leq \gn^*.$ Thus
$\sup\{\gm_1,\gm_2\}\leq \gn^*$ which reads $\gn\leq\gn^*$, and
equality follows. Next, assume $\gth\in[0,1]$. Then
$\gm_\gth=\gth\gm_1+(1-\gth)\gm_2\leq\gn=\sup\{\gm_1,\gm_2\}$.
Since $\gn\in\CG^\Gw(g)$, and any measure dominated by a good
measure is a good measure, $\gm_\gth\in\CG^\Gw(g)$. It follows by
\rth{redmes2} that $\gm_\gth=\gm_\gth^*$.\smallskip

\noindent Next, by \rcor{redmes4'}, $[\inf\{\gm^*,\gn^*\}]$ is a
good measure. Since $[\inf\{\gm^*,\gn^*\}]\leq [\inf\{\gm,\gn\}]$,
it follow by \rth{redmes2} that
\begin{equation}\label{inf1}
\inf\{\gm^*,\gn^*\}\leq [\inf\{\gm,\gn\}]^*.
\end{equation}
Conversely,
$$\inf\{\gm,\gn\}\leq \gm\Longrightarrow[\inf\{\gm,\gn\}]^*\leq\gm^*,
$$
and similarly with $\gn$. Thus $ [\inf\{\gm,\gn\}]^*\leq \inf\{\gm^*,\gn^*\}$. \smallskip

\noindent For the last assertion, by Hahn's decomposition theorem there exist two disjoint Borel sets $A$ and $B$ such that $\Gw=A\cup B$ and
$\sup\{\gm,\gn\}=\chi_{_{A}}\gm+\chi_{_{B}}\gn$. Actually, $\gm\geq\gn$ on $A$ and $\gn\geq\gm$ on $B$. This implies also $\sup\{\gm^*,\gn^*\}=\chi_{_{A}}\gm^*+\chi_{_{B}}\gn^*$. Thus, by \rcor {redmes4'},
$$[\sup\{\gm,\gn\}]^*=(\chi_{_{A}}\gm+\chi_{_{B}}\gn)^*=\chi_{_{A}}\gm^*+\chi_{_{B}}\gn^*=\sup\{\gm^*,\gn^*\},$$
since $\sup\{\chi_{_{A}}\gm^*,\chi_{_{B}}\gn^*\}=\sup\{\gm^*,\gn^*\}$.
\qeda
\bth{stability} Let $\gm,\gn\in\frak M_{+}^1$. Then
\begin {equation}\label{contract}
\abs {\gm^*-\gn^*}\leq \abs {\gm-\gn}.
\end {equation}
\es
\Proof We first assume $\gm\geq\gn$. By \rth{ord},
$$0\leq \gm^*-\gn^*\leq\gm-\gn.
$$
This implies (\ref{contract}). Next we write $\sup\{\gm,\gn\}=\gn+(\gm-\gn)_{+}$. Since $\gn\leq
\sup\{\gm,\gn\}$, $\gn^*\leq [\sup\{\gm,\gn\}]^*=\sup\{\gm^*,\gn^*\}$ by \rth{structh}. Thus
$$[\sup\{\gm,\gn\}]^*-\gn^*\leq \sup\{\gm,\gn\}-\gn=(\gm-\gn)_{+}.
$$
Thus implies $(\gm^*-\gn^*)_{+}\leq (\gm-\gn)_{+} $. Similarly $(\gn^*-\gm^*)_{+}\leq (\gn-\gm)_{+} $.
\qeda\\

In order to characterize the universally good measures, we introduce a capacity naturally associated to the weak formulation of problem (\ref{Eq2-1}). This yields to a capacity type characterization of $H^N$.
If $K\subset\Gw$ is compact, we denote
\begin {equation}\label{equa0}\BA {l}
c_\Gw(K)=\inf\left\{\dint_{Q_T}\abs{\prt_t\psi+\Gd\psi}\,dx\,dt:\right.\\
\phantom{\myfrac{A}{B}----------}
\left.\psi\in
C^{2,1}_{\ell,0}(\bar Q_T),\,\psi(x,0)\geq 1\txt {in a neighborhood of}K\phantom{\myfrac{A}{B}}\!\!\!\!\!\!\right\}.
\EA\end {equation}
\bth {equa1} For every compact $K\subset\Gw$, we have
\begin {equation}\label{equiv1}
H^N(K)=c_\Gw(K).
\end {equation}
\es
\Proof Let $K\subset\Gw$ be compact.\medskip

\noindent{\it Step 1. }We claim that for any $\ge>0$, there exists
$\psi_\ge=\psi\in C^{2,1}_{\ell,0}(\bar Q_T)$ such that $\psi\geq
0$ in $Q_T$, $\psi(x,0)\geq 1\txt {on}K$ and
\begin{equation}\label{equa2}
\dint_{Q_T}\abs{\prt_t\psi+\Gd\psi}\,dx\,dt\leq c_\Gw(K)+\ge.
\end{equation}
Let $\xi\in C^{2,1}_{\ell,0}(\bar Q_T)$ such that $\xi(x,0)\geq 1\txt {on}K$ and
$$
\dint_{Q_T}\abs{\prt_t\xi+\Gd\xi}\,dx\,dt\leq c_\Gw(K)+\ge/2.
$$
Let $\{\eta_j\}$ be a regularizing sequence depending only on the space variable and such that the support of $\eta_j$ is contained in the ball $B_{\ge_j}$, with $\ge_j\to 0$ as $j\to\infty$. If we extend $\xi$ in $\BBR^N\ti [0,T]$ as a $C^{2,1}$-function, we set
$$f_j(x,t)=\eta_j\ast\abs{\prt_t\xi+\Gd\xi}(x,t)=\int_{\Gw}\eta_j(x-y)\abs{\prt_t\xi+\Gd\xi}(y,t)dy.
$$
If $j\to\infty$, $\{f_j\}$ converges to $\abs{\prt_t\xi+\Gd\xi}$ uniformly in $\bar Q_T$. Let $v_j$ be the solution of
$$\left\{\BA {l}
\prt_tv_j+\Gd v_j=-f_j\quad \text {in } Q_T\\[2mm]
\phantom{\prt_tv_j+\Gd }
v_j=0\quad \text {in } \prt_\ell Q_T\\[2mm]
\phantom{\prt_tv_j}
v_j(.,T)=0\quad \text {in } \Gw.
\EA\right.
$$
Clearly $v_j\geq 0$ in $Q_T$. Let $v$ be the solution of
$$\left\{\BA {l}
\prt_tv+\Gd v=-\abs{\prt_t\xi+\Gd\xi}\quad \text {in } Q_T\\[2mm]
\phantom{\prt_tv+\Gd }
v=0\quad \text {in } \prt_\ell Q_T\\[2mm]
\phantom{\prt_tv}
v(.,T)=0\quad \text {in } \Gw.
\EA\right.
$$
By the maximum principle $v\geq \max\{\xi,0\}$, thus
$v(x,0)\geq 1$ on $K$. Because $v_j(x,0)\to v(x,0)$ uniformly on $\bar\Gw$, for any $0<\ga<1$, we can fix $j_\ga$ such that $v_{j_\ga}(x,0)\geq \ga$ on $K$ and
$\norm{f_{j_\ga}}_{L^1(Q_T)}\leq \norm{\prt_t\xi+\Gd\xi}_{L^1(Q_T)}+\ge/4$. Next
$\psi_\ga=\ga^{-1} v_{j_\ga}$. Then $\psi_\ga\geq 0$ in $Q_T$, and $\psi_\ga (x,0)\geq 1$ on $K$. Moreover
$$\BA {l}\dint_{Q_T}\abs{\prt_t\psi_\ga+\Gd\psi_\ga}\,dx\,dt
=\ga^{-1}\dint_{Q_T}\abs{\prt_tv_{j_\ga}+\Gd v_{j_\ga}}\,dx\,dt\\[2mm]
\phantom{\dint_{Q_T}\abs{\prt_t\psi_\ga+\Gd\psi_\ga}\,dx\,dt}
\leq \ga^{-1}\left(\dint_{Q_T}\abs{\prt_t\xi+\Gd\xi}\,dx\,dt+\ge/4\right)\\[2mm]
\phantom{\dint_{Q_T}\abs{\prt_t\psi_\ga+\Gd\psi_\ga}\,dx\,dt}
\leq \ga^{-1}\left(c_\Gw(K)+3\ge/4\right).
\EA$$
Next we fix
$$\ga=\myfrac{c_\Gw(K)+3\ge/4}{c_\Gw(K)+\ge}$$
and derive (\ref{equa2}).\smallskip

\noindent{\it Step 2. } There holds
\begin{equation}\label {ineq1}
H^N(K)\leq c_\Gw(K).
\end {equation}
From (\ref{equa2}),
$$\dint_{Q_T}\left(-\prt_t\psi-\Gd\psi\right)\,dx\,dt\leq \dint_{Q_T}\abs{\prt_t\psi+\Gd\psi}\,dx\,dt\leq c_\Gw(K)+\ge.
$$
But
$$\BA {l}\dint_{Q_T}\left(-\prt_t\psi-\Gd\psi\right)\,dx\,dt=
\myint{\Gw}{}\psi(x,0)dx-\dint_{\prt_\ell Q_T}\myfrac{\prt\psi}{\prt n}dS\,dt
\geq H^N(K)
\EA$$
since $\psi(x,T)=0$, $\psi(x,0)\geq 1$ on $K$, and the normal derivative of $\psi$ on
$\prt_\ell Q_T$ is nonpositive. This yields to (\ref{ineq1}) because $\ge$ is arbitrary.
\smallskip

\noindent{\it Step 3. }For any $\ge >0$ there exists $\psi\in C^{2,1}_{\ell,0}(\bar Q_T)$ such that $0\leq \psi\leq 1+\ge$ in $Q_T$, $\psi (x,0)\geq 1$ on $K$ and
\begin{equation}\label{equa3}
\dint_{Q_T}\abs{\prt_t\psi+\Gd\psi}\,dx\,dt\leq H^N(K)+\ge.
\end{equation}
For $\gd>0$ let $K_\gd=\{x\in\BBR^N:\dist(x,K)\leq\gd\}$. By the regularity of $H^N$, we can choose $\gd$ small enough such that
$$H^N(K_\gd\cap\Gw)\leq H^N(K)+\ge/5.
$$
We fix $\xi\in C^{2,1}_{\ell,0}(\bar Q_T)$ such that $0\leq\xi\leq 1$ and
$$\xi(x,0)=\left\{\BA {l}1\;\text { if }x\in K_{\gd/2}\\[2mm]
0\;\text { if }x\in \bar\Gw\setminus K_{\gd}.
\EA\right.$$
Let $\gs>0$ and
$$\gr_\gs(x,t)=\left(1-\myfrac{t}{\gs}\right)_+^2.
$$
Since $\abs{\xi_t+\Gd\xi}(x,t)=0$ a. e. on $\{(x,t):\xi(x,t)=0\}$ and
$\{(x,t):\gr_\gs(x,t)>0\}\subset \bar\Gw\ti [0,\gs]$,
we can choose $\gs$ such that
$$\dint_{\prt_\ell Q_T\cap \{(x,t):\xi\leq\gr_\gs\}}\myfrac{\prt\xi}{\prt n}dS\,dt
+\dint_{ \{(x,t):\xi\leq\gr_\gs\}}\abs{\xi_t+\Gd\xi}\,dx\,dt\leq\ge/5.
$$
We set $u=\gr_\gs-(\gr_\gs-\xi)_+$. Because $\gr_\gs$ is independent of $x$, the argument developed by Brezis and Ponce \cite{BP} applies in the sense that $\Gd u(.,t)\in\frak M(\Gw)$ and $\Gd u(.,t)=\Gd \xi(.,t)$ on $\{x:\xi(x,t)<\gr_\gs(t)\}$ and more explicitely
$\prt_tu+\Gd u=\prt_t\xi+\Gd\xi$ on $\{(x,t):\xi(x,t)<\gr_\gs(t)\}$. In addition
$$
\prt_tu=\prt_t\gr_\gs-\text{sign}_+(\gr_\gs-\xi)(\prt_t\gr_\gs-\prt_t\xi),
$$
and $\prt_tu=\prt_t\gr_\gs$ a.e. on $\{(x,t):\xi(x,t)\geq \gr_\gs(x,t)\}$. Because $\gr_\gs$ is decreasing, we finally obtain
$$\prt_tu+\Gd u\leq 0\text { on }\{(x,t):\xi(x,t)\geq \gr_\gs(x,t)\}.
$$
We notice that $\prt_tu$ is bounded, and, following \cite {BP},
\begin{equation}\label{meas1}\BA {l}
\norm{\prt_tu+\Gd u}_{\frak M}=
\norm{(\prt_tu+\Gd u)\chi_{_\{\xi\geq\gr_\gs\}}}_{\frak M}+
\norm{(\prt_tu+\Gd u)\chi_{_\{\xi<\gr_\gs\}}}_{\frak M}\\[2mm]
\phantom{\norm{\prt_tu+\Gd u}_{\frak M}}
=\norm{(\prt_tu+\Gd u)\chi_{_\{\xi\geq\gr_\gs\}}}_{\frak M}
+\dint_{\{\xi<\gr_\gs\}}\abs {\prt\xi+\Gd\xi}\,dx\,dt\\[2mm]
\phantom{\norm{\prt_tu+\Gd u}_{\frak M}}
=-\dint_{\{\xi\geq\gr_\gs\}}d(\prt_tu+\Gd u)+\dint_{\{\xi<\gr_\gs\}}\abs {\prt\xi+\Gd\xi}\,dx\,dt\\[2mm]
\phantom{\norm{\prt_tu+\Gd u}_{\frak M}}
\leq
-\dint_{Q_T}d(\prt_tu+\Gd u)+2\dint_{\{\xi<\gr_\gs\}}\abs {\prt\xi+\Gd\xi}\,dx\,dt\\[2mm]
\phantom{\norm{\prt_tu+\Gd u}_{\frak M}}
\leq -\dint_{Q_T}d(\prt_tu+\Gd u)+2\ge/5.
\EA\end{equation}
Next, by definition,
\begin{equation}\label{meas2}\BA {l}
-\dint_{Q_T}d(\prt_tu+\Gd u)=
-\dint_{Q_T}u(\prt_t1-\Gd1)\,dx\,dt-\myint{\Gw}{}(u(x,T)-u(x,0))dx          \\[2mm]
\phantom{-\dint_{Q_T}d(\prt_tu+\Gd u)=
-\dint_{Q_T}u(\prt_t1-\Gd1)\,dx\,dt}
-\dint_{\prt_\ell Q_T}\myfrac{\prt u}{\prt n}dS\,dt
\\[2mm]
\phantom{-\dint_{Q_T}d(\prt_tu+\Gd u)}
=\myint{\Gw}{}u(x,0)dx-\dint_{\prt_\ell Q_T}\myfrac{\prt u}{\prt n}dS\,dt
\\[2mm]
\phantom{-\dint_{Q_T}d(\prt_tu+\Gd u)}
=\myint{\Gw}{}\xi(x,0)dx+\dint_{\prt_\ell Q_T\cap \{(x,t):\xi\leq\gr_\gs\}}\myfrac{\prt\xi}{\prt n}dS\,dt\\[2mm]
\phantom{-\dint_{Q_T}d(\prt_tu+\Gd u)}
\leq H^N(K_\gd)/+\ge/5\\[2mm]
\phantom{-\dint_{Q_T}d(\prt_tu+\Gd u)}
\leq H^N(K)+2\ge/5.
\EA\end{equation}
We finally derive
\begin{equation}\label{meas3}
\norm{\prt_tu+\Gd u}_{\frak M}\leq H^N(K)+4\ge/5.
\end{equation}
Next, we smooth the measure $\abs{\prt_tu+\Gd u}$ using a space convolution process with the same $\eta_{j}$, as in Step 1.  One can construct a function
$\psi\in C^{2,1}_{\ell,0}(\bar Q_T)$ such that $ 0\leq \psi\leq 1+\ge$ in $\bar Q_T$,
$\psi (x,0)\geq 1$ on $K$ and
\begin{equation}\label{meas4}
\dint_{Q_T}\abs{\prt_t\psi+\Gd \psi}dx\,dt\leq \norm{\prt_tu+\Gd u}_{\frak M}+\ge/5.
\end{equation}
Combining (\ref{meas3}) and (\ref{meas4}), one derive (\ref{equa3}).\smallskip

\noindent{\it Step 4. }There holds
 \begin{equation}\label{ineq2}
c_\Gw(K)\leq H^N(K).
\end{equation}
Actually, (\ref{equa3}) implies
$$c_\Gw(K)\leq H^N(K)+\ge.
$$
Letting $\ge\to 0$ yields to (\ref{ineq2}).\qeda\\

Thanks to this result we are able to characterize the universally good measures.
\bth{uni} Let $\gm\in \frak M_+^1(\Gw)$. If $\gm\in\CG^\Gw(g)$ for any function $g$ satisfying (\ref{g1}), then $\gm\in L^1_\gr(\Gw)$.
\es
\Proof We follow essentially the proof of \cite[Th 7]{BP}.\medskip

\noindent {\it Step 1.} We claim that for every Borel set $\Gs\subset\Gw$, such that $H^N(\Gs)=0$, there exists a continuous function $g$ verifying (\ref{g1}) such that $\gm^*=0$ for any $\gm\in\frak M^1_+(\Gw)$ satisfying $\gm(\Gs^c)=0$.\medskip

\noindent Let $\{K_{j}\}_{j\in\BBN^*}$ be an increasing sequence
of compact subsets of $\Gs$ such that $K=\cup_{j}K_{j}$ and
$\gm(\Gs\setminus K)=0$. Since $H^N(K_{j})=0$ for any $j\geq 1$,
it follows from \rth{equa1} [Step 3], that there exists
$\psi_{j}\in C^{2,1}_{{\ell,0}}(\bar Q_{T})$ such that $0\leq
\psi_{j}\leq 2$ in $Q_{T}$, $\psi_j(x,0)\geq 1$ on $K_{j}$ and
$$\dint_{Q_{T}}\abs{\prt_{t}\psi_{j}+\Gd\psi_{j}}dx\,dt\leq 1/j.
$$
In particular,
$$\abs{\prt_{t}\psi_{j}+\Gd\psi_{j}}\to 0\quad \text{a.e. in }Q_{T},
$$
and, since $\psi_{j}$ solves
$$\left\{\BA{l}\prt_{t}\psi_{j}+\Gd\psi_{j} =\ge_{j}\quad\text {in }Q_{T}\\[2mm]
\phantom{\prt_{t}\psi_{}}
\psi_{j}    (x,T)=0\quad\text {in }\Gw\\[2mm]
\phantom{\prt_{t}\psi_{,}}
\psi_{j}    (x,t)=0\quad\text {on }\prt_{\ell}Q_{T}
\EA\right.$$
with $\ge_{j}\to 0$ in $L^1(Q_{T})$ it follows $\psi_{j}\to 0$ in $L^1(Q_{T})$ and a.e.. Furthermore there exists some $G\in L_{\gr}^1(Q_{T})$ such that
$$\gr^{-1}\abs{\prt_{t}\psi_{j}+\Gd\psi_{j}}\leq G\forevery j\in\BBN^*.
$$

By a theorem of De La Vall\'ee-Poussin noticed in \cite{P}, there
exists a convex function $h:(-\infty,\infty)\mapsto [0,\infty)$
such that $h(s)=0$ for $s\leq 0$, $h(s)>0$ for $s>0$,
$$\lim_{t\to\infty}\myfrac{h(t)}{t}=\infty\quad\text{and }h(G)\in L_{\gr}^1(Q_{T}).
$$
Let $g=h^*$ be the convex conjugate of $h$. We denote by $\gm^*=\gm^*(g)$ the reduced measured associated to $g$. Since $\gm^*\in \CG^\Gw(g)$, we denote by $u$ the solution of the corresponding initial value problem. Taking $\psi_{j}$ as a test function in (\ref {Eq2-1}), we obtain
\begin{equation}\label{uni1}
\dint_{Q_{T}}\left(-u(\prt_{t}\psi_{j}+\Gd\psi_{j})+\psi_{j}g(u)\right)dx\,dt=\myint{\Gw}{}\psi_{j}(x,0)d\gm^*.
\end{equation}

We first assume that $\gm\in\frak M^0(\Gw)$, thus we can take $1$ as a test function (this is easily justified by approximations) and obtain
\begin{equation}\label{uni2}
\dint_{Q_{T}}g(u)dx\,dt=\myint{\Gw}{}d\gm^*.
\end{equation}
Therefore
\begin{equation}\label{uni3}
\gm^*(K_{j})\leq \dint_{Q_{T}}\left(-u(\prt_{t}\psi_{j}+\Gd\psi_{j})+\psi_{j}g(u)\right)dx\,dt
\end{equation}
and
\begin{equation}\label{uni4}\BA {l}
\abs{-u(\prt_{t}\psi_{j}+\Gd\psi_{j})+\psi_{j}g(u)}
\leq \myfrac{\abs{\prt_{t}\psi_{j}+\Gd\psi_{j}}}{\gr}u\gr+\psi_{j}g(u)\\[2mm]
\phantom{\abs{-u(\prt_{t}\psi_{j}+\Gd\psi_{j})+\psi_{j}g(u)}}
\leq h\left(\gr^{-1}\abs{\prt_{t}\psi_{j}+\Gd\psi_{j}}\right)\gr+g(u)\gr+\psi_{j}g(u)
\\[2mm]
\phantom{\abs{-u(\prt_{t}\psi_{j}+\Gd\psi_{j})+\psi_{j}g(u)}}
\leq h(G)\gr+ Cg(u)
\EA\end{equation}
By Lebesgue's theorem, the right-hand side of (\ref{uni3}) tends to $0$ when $j\to\infty$. Thus $\gm^*(K_{j})=0$, for any $j\in\BBN^*$, and finally $\gm^*(\Gs)=0$.

Next we assume $\gm\in \frak M^1(\Gw)$. Then there exists an increasing sequence of $\gm_{n}\in \frak M^0(\Gw)$ with compact support in $\Gw$ such that $\gm_{n}\uparrow\gm$. Using what is proved above, $\gm_{n}^*(\Gs)=0$ and, by \rth{ord}, $\gm^*\leq \gm-\gm_{n}$, thus
$\gm^*(\Gs)\leq (\gm-\gm_{n})(\Gs)$. Letting $n\to\infty$ implies $\gm^*(\Gs)=0$.\smallskip

\noindent {\it Step 2.} If $\gm\in\frak M^1_{+}(\Gw)$ is good, for any Borel set $\Gs\subset\Gw$, with $H^{N-1}(\Gs)=0$, we denote $\gn=\gm\vline_{\Gs}$. Then there exists $g_{\gn}$ such that $g^*_{\gn}=0$. Since $\gn\leq\gm$, $\gn\in\CG^\Gw(g)$, thus $\gn=\gn^*=0$ and finally, $\gm(\Gs)=0$. Thus
$\gm\in L^1_{\gr}(\Gw)$.\qeda
\section{The Cauchy-Dirichlet problem}
\setcounter{equation}{0}
In this section $\Gw$ is again a smooth bounded domain in $\BBR^N$ and $\gr(x)=\dist (x,\prt\Gw)$. We denote by $\mathfrak M(\prt_{\ell}Q_{T})$  the set of Radon measures in $\prt_{\ell}Q_{T}$ and by $\mathfrak M_{+}(\prt_{\ell}Q_{T})$, the positive ones. The function $g$ is supposed to satisfy (\ref{g1}). We consider the Cauchy-Dirichlet problem
\begin {equation}\label {CD1}\left\{\BA {l}
\prt_tu-\Gd u+g(u)=0\quad \text {in } Q_T:=\Gw\ti(0,T)\\[2mm]
\phantom{\prt_t-\Gd u+g(u)}
u=\gm\quad \text {in } \prt_\ell Q_T:=\prt\Gw\ti(0,T)\\[2mm]
\phantom{\Gd u+g(u)}
u(.,0)=0\quad \text {in } \Gw,
\EA\right.
\end {equation}

\bdef {solCD} Let $\gm\in\mathfrak M_{+}(\prt_{\ell}Q_{T})$. A function $u\in L^1(Q_{T})$ is a weak solution of (\ref{CD1}) if $g(u)\in L^1_{\gr}(Q_{T})$ and
\begin{equation}\label {CD2}
\dint_{Q_T}\left(-u\prt_t\gz-u\Gd\gz+\gz\,g(u) \right)dx=-\myint{\prt_{\ell}Q_{T}}{}\myfrac{\prt\gz}{\prt\gn}d\gm,
\end{equation}
for every $\gz\in C^{2,1}_{0}(\bar Q_{T})$. \es Solutions of
(\ref{CD1}) are always unique; sufficient conditions for existence
are developed in \cite{MV2}. We define, similarly to the cases of
the initial value problem, super and subsolutions of \ref{CD1}. In
which case, the equality sign in \ref {CD2} is replaced by $\geq$
and $\leq$ respectively, the integrability conditions on $u$ and
$g(u)$ being preserved. As simple example for existence of a
solution it is the case when $g$ satisfies
\begin{equation}\label {CD3}
\dint_{Q_T}g(\BBP^H[\gm](x,t))\gr(x)dx\,dt<\infty.
\end{equation}
In this formula $\BBP^H[\gm]$ is the Poisson-heat potential of $\gm$ in $Q_T$, that is the solution of
\begin{equation}\label {CD4}\left\{\BA {l}
\prt_{t}v-\Gd v=0\qquad\text{in }Q_{T}\\[2mm]
\phantom{\prt_{t}v-\Gd }v=\gm\qquad\text{in }\prt_{\ell}Q_{T}\\[2mm]
\phantom{\prt_{t}v-\Gd }v=0\qquad\text{in }\Gw.
\EA\right.\end{equation}

\bdef{good}A measure $\gm$ for which problem (\ref{CD1}) can be solved is called a good measure relative to $g$ for the Cauchy-Dirichlet problem. The set of good measures is denoted by
$\CG^{\prt_{\ell}Q_{T}}(g)$, and a universally good measure is a measure which belongs to
$\CG^{\prt_{\ell}Q_{T}}(g)$ for any $g$ satisfying (\ref{g1}).
\es

The notion of lateral trace is defined in \cite{MV3}. For $\gb>0$,
we denote
$$\Gw_{\gb}=\{x\in\Gw:\gr(x)<\gb\},\;\Gw'_{\gb}=\{x\in\Gw:\gr(x)>\gb\}\text { and }\Gs_{\gb}=\prt\Gw_{\gb}.
$$
We shall also denote $\Gs=\Gs_{0}=\prt\Gw$. There exists $\gb_{0}>0$ such that for any $\gb\in (0,\gb_{0}]$, the mapping
$x\in \Gw_{\gb}\mapsto (\gs(x),\gr(x)$, where $\gs(x)$ is the unique point on $\prt\Gw$ which minimizes
the distance from $x$ to $\prt\Gw$, is a $C^2$ diffeomorphism from $\bar\Gw_{\gb}$ to $\Gs\ti [0,\gb_{0}]$. If $\gf\in L^1_{loc}(\prt_{\ell}Q_{T})$, we denote $\gf^\gb(x,t)=\gf (\gs(x),t)$, for any $x\in\Gs_{\gb}$ and $dS_\gb$ is the surface measure on $\Gs_{\gb}$.
for the sake of simplicity

\bdef{lat-trace} Let $u\in L^1_{loc}(Q_{T})$. 1- We say that $u$ admits the Radon measure $\gm\in\frak M_{+}(\prt_{\ell}Q_{T})$ as a lateral boundary trace if  it exists
\begin{equation}\label {tr-lat}
\rm ess\,lim_{\gb\to 0}\int_{0}^T\int_{\Gs_{\gb}}u\gf^\gb\,dS_{\gb}dt=\dint_{\prt_{\ell}Q_{T}}\gf\,d\gm
\forevery\gf\in C_{0}(\CR).
\end {equation}
We shall denote $\gm=Tr_{\prt_{\ell}Q_{T}}(u)$.\smallskip

\noindent 2- We say that $u$ admits the outer regular Borel measure
$\gn\approx (\Gs,\gm)$ as a lateral boundary trace if  it exists an open subset $\CR\subset \prt_{\ell}Q_{T}$ and $\gm\in \frak M_{+}(\CR)$ such that
\begin{equation}\label {tr-lat''}
\rm ess\,lim_{\gb\to 0}\int_{0}^T\int_{\Gs_{\gb}}u\gf^\gb\,dS_{\gb}dt=\infty
\forevery\gf\in C_{0}(\prt_{\ell}Q_{T}),\gf\geq 0,\gf>0\text{ somewhere in } \CS,
\end {equation}
with $\CS=\prt_{\ell}\setminus\CR$.
We shall denote $\gn=Tr_{\prt_{\ell}Q_{T}}(u)$.
\es
Propositions 2.5, 2.6, 2.7, 2.8 and \rth{trsub} are still valid, if we replace the notion of initial trace by the notion of lateral boundary trace.The new version of \rth{trsub} is the following.

\bth {Ltrsub} Let $u$ be a nonnegative subsolution of (\ref{CD1}). Then the lateral boundary trace of $u$
is a positive Radon measure $\tilde \gm$ such that $\tilde\gm\leq \gm$. Furthermore, if (\ref{CD1}) admits a weak solution $u_{\gm}$ there holds $u\leq u_{\gm}$.
\es

We consider now a sequence of functions $g_{k}$ satisfying (\ref{Eq2-4}). For any positive Radon measure $\gm$ on $\prt_{\ell}Q_{T}$, it is possible to solve, with $u=u_{k}$,
\begin {equation}\label {LEq2-5}\left\{\BA {l}
\prt_tu-\Gd u+g_k(u)=0\quad \text {in } Q_T\\[2mm]
\phantom{\prt_t-\Gd u+g_k(u)}
u=\gm\quad \text {in } \prt_\ell Q_T\\[2mm]
\phantom{\Gd u+g_k(u)}
u(.,0)=0\quad \text {in } \Gw.
\EA\right.
\end {equation}

The following result is proved as \rth{redmes}

\bth {Lredmes} When $k\to\infty$, the sequence $\{u_k\}$ converges in $L^1(Q_T)$ to a some nonnegative function $u^*$ such that
$g(u^*)\in L_\gr^1(Q_T)$ and there exists a positive Radon measure $\gm^*$ smaller that $\gm$ with the property that
\begin {equation}\label {LEq2-6}\left\{\BA {l}
\prt_tu^*-\Gd u^*+g(u^*)=0\quad \text {in } Q_T\\[2mm]
\phantom{\prt_t-\Gd u^*+g(u^*)}
u^*=\gm^*\quad \text {in } \prt_\ell Q_T\\[2mm]
\phantom{\Gd u^*+g(u^*)}
u^*(.,0)=0\quad \text {in } \Gw.
\EA\right.
\end {equation}
Furthermore $u^*$ is the largest subsolution of problem (\ref{CD1}).
\es

Mutadis mutandis, the reduced measure $\gm^*$ on the lateral boundary inherits the properties of the reduced measure at initial time and the assertions of Theorems 2.13, 2.14, Corollaries 2.15, 2.16, \rth {ord}, Corollaries 2.18, 2.19 and Theorems 2.20 and 2.21, are valid in the framework of the lateral bondary reduced measure. The main novelty is the intruction of a new capacity on $\prt_{\ell}Q_{T}$. If $K\subset\prt_{\ell}Q_{T}$ is compact, we denote
\begin {equation}\label{Lequa0}\BA {l}
c_{_{\prt_{\ell}Q_{T}}}(K)=\inf\left\{\dint_{Q_T}\abs{\prt_t\psi+\Gd\psi}\,dx\,dt\right.\\[2mm]
\phantom{----------}\left.:\psi\in
C^{2,1}_{\ell,0}(\bar Q_T),\,-\myfrac{\prt\psi}{\prt\gn}\geq 1\txt {in some neighborhood of}K\phantom{\myfrac{A}{B}}\!\!\!\!\!\!\right\}.
\EA\end {equation}

\bth {Lequa1} For every compact $K\subset\prt_{\ell}Q_{T}$, we have
\begin {equation}\label{equiv1*}
H^N(K)=c_{_{\prt_{\ell}Q_{T}}}(K).
\end {equation}
\es

\Proof Let $K\subset\prt_{\ell}Q_{T}$ be compact.\smallskip

\noindent {\it Step 1}. For any $\ge>0$ there exists $\psi\in C^{2,1}_{\ell,0}(\bar Q_T)$, $\psi\geq 0$ such that $-\prt\psi(x,t)/\prt\gn\geq 1$ in some neighborhood of $K$.\smallskip

\noindent Let $\xi\in C^{2,1}_{\ell,0}(\bar Q_T)$ such that $-\prt\psi(x,t)/\prt\gn\geq 1\txt {on}K$ and
$$\dint_{Q_T}\abs{\prt_t\xi+\Gd\xi}\,dx\,dt\leq c_{_{\prt_{\ell}Q_{T}}}(K)+\ge/2.
$$
We extend $\xi$ as a $C^{2,1}(\BBR^N\ti [0,T])$-function and define $f_{j}$, $v_{j}$ and $v$ in the same way as in the proof of \rth {equa1}, Step 1. Since $f_{j}\to \prt_t\xi+\Gd\xi$ uniformly in $\bar Q_{T}$,
$$\myfrac{\prt v_{j}}{\prt\gn}\to \myfrac{\prt v}{\prt\gn},
$$
uniformly in $\bar Q_{T}$. Since $v$ and $\xi$ vanishes on $\prt_{\ell}Q_{T}$ and at $t=T$, $v\geq\xi$, thus
$$0\leq \myfrac{\prt \xi}{\prt\gn}\leq -\myfrac{\prt v}{\prt\gn}\quad\text{ on }\prt_{\ell}Q_{T},
$$
and $-\prt v/\prt\gn\geq1$ in some neighborhood of $K$. For $\ga\in (0,1)$ we fix $j_{0}$ such that
$-\prt v_{j_{0}}/\prt\gn\geq\ga$ on $K$ and
$$\dint_{Q_T}\abs{\prt_tv_{j_{0}}+\Gd v_{j_{0}}}\,dx\,dt\leq \dint_{Q_T}\abs{\prt_t\xi+\Gd\xi}\,dx\,dt+ \ge/4.
$$
We set $\psi=\ga^{-1} v_{j_{0}}$ and get
$$\dint_{Q_T}\abs{\prt_t\xi+\Gd\xi}\,dx\,dt\leq \ga^{-1}(c_{_{\prt_{\ell}Q_{T}}}(K)+3\ge/4).
$$
We end the proof as in \rth {equa1}, Step 1.\smallskip

\noindent {\it Step 2}. In this step we follow essentially the proof of \cite[Lemma 8]{BP}. For any $\ge>0$ there exists $\psi\in C^{2,1}_{\ell,0}(\bar Q_T)$, such that $0\leq \psi\leq \ge$, $-\prt\psi(x,t)/\prt\gn\geq 1$ in some neighborhood of $K$ and
\begin{equation}\label{B1}\dint_{Q_T}\abs{\prt_t\psi+\Gd\psi}\,dx\,dt\leq H^N(K)+\ge\quad\text {and }\abs{\myfrac{\psi}{\gr}}\leq 1+\ge\text { in }Q_{T}.
\end {equation}
Let $\gd>0$ and $\tilde N_\gd(K)=\{(x,t):\dist ((x,t),K)\}$, be such that
$$H^{N}(N_{\gd}(K)\cap\prt_{\ell}Q_{T})\leq H^{N}(K)+\ge
$$
We take $\xi\in C^{2,1}_{\ell,0}(\bar Q_T)$ such that $\xi>0$ in $Q_T$, $\prt\xi/\prt\gn=-1$ on $N_{\gd/2}(K)\cap\prt_\ell Q_T$ and $\prt\xi/\prt\gn=0$ on $\prt_\ell Q_T\setminus N_{\gd}$, $0\leq-\prt\xi/\prt\gn\leq 1$ and $\xi/\gr\leq 1+\ge$, we first take $a>0$ small enough so that
$$\dint_{\prt_\ell Q_T\cap\{\xi<a\}}\myfrac{\prt\xi}{\prt\gn}dS\,dt
+\dint_{Q_T\cap\{\xi<a\}}\abs{\prt_t\xi+\Gd\xi}dx\,dt<\ge,
$$
and set $u=a-(a-\gz)_+$. Then, the same method as in \rth {equa1}-Step 3 yields to
\begin {equation}\label {Lmeas3}
\norm{\prt_tu+\Gd u}_{\mathfrak M}\leq H^N(K)+4\ge/5.
\end {equation}
The conclusion of the proof is similar.\qeda\\

By an easy adaptation of the proof of \rth{uni} we have the following characterization of the universally good measures.

\bth{Luni} Let $\gm\in \frak M_+(\prt_\ell Q_T)$. If $\gm\in\CG^{\prt_\ell Q_T}(g)$ for any function $g$ satisfying (\ref{g1}), then $\gm\in L^1(\prt_\ell Q_T)$.
\es

\end{document}